\documentclass[11pt]{article}

\usepackage{fullpage,graphicx,authblk,ifthen}

\usepackage{titlesec}

\usepackage{amsthm}
\usepackage{mathtools}
\usepackage{accents}	
\newtheorem{theorem}{Theorem}
\newcounter{def}
\newtheorem{definition}[def]{Definition}
\newcounter{examples}
\newenvironment{example}[1][Example \theexamples]
               {\refstepcounter{examples}
                 \begin{trivlist} \item[\hskip \labelsep {\bfseries #1}]}
               {\qed \end{trivlist}}
               
\usepackage[authoryear]{natbib}
\usepackage{amssymb}
\usepackage{accents}
\usepackage[dvipsnames,table,xcdraw]{xcolor}		
\usepackage{todonotes}	

\usepackage[funcfont=italic,full]{./complexity/complexity}	

\usepackage{hyperref}

\usepackage{subcaption}	

\newcommand{\proj}{\operatorname{proj}}

\newcommand{\Z}{\mathbb{Z}}
\newcommand{\Q}{\mathbb{Q}}
\newcommand{\B}{\mathbb{B}}

\renewcommand{\Re}{\mathbb{R}}

\newcommand{\Dmcal}{\mathcal{D}}
\newcommand{\Fmcal}{\mathcal{F}}
\newcommand{\Emcal}{\mathcal{E}}

\newcommand{\phiSS}{f_x}
\newcommand{\uphiSS}{\underaccent{\bar}{f}_x}
\newcommand{\phiVF}{\phi}
\newcommand{\uphiVF}{\underline{\phi}}
\newcommand{\uphi}{\underline{\phi}}
\newcommand{\phiIP}{\phi_{IP}}
\newcommand{\uphiIP}{\underline{\phi}_{IP}}
\newcommand{\bphiIP}{\bar{\phi}_{IP}}
\newcommand{\bphiNK}{\bar{\phi}_{N \setminus K}}
\newcommand{\bphi}[1]{\bar{\phi}^{#1}_{N \setminus K}}
\newcommand{\urho}{\underline{\rho}}

\newcommand{\Pcal}{{\cal P}}

\newcommand{\ra}{{\rightarrow}}

\newcommand{\sm}{{\setminus}}

\newcommand{\stkout}[1]{\ifmmode\text{\sout{\ensuremath{#1}}}\else\sout{#1}\fi}
\newcommand{\utilde}[1]{\underaccent{\tilde}{#1}}

\newcommand{\noprint}[1]{}
\newcommand{\midd}{\;\middle|\;}
\renewcommand{\P}{\mathcal{P}}

\newtheorem{assumption}{Assumption}

\begin{document}

\title{A Framework for Generalized Benders' Decomposition and Its
  Application to Multilevel Optimization}

\author{Suresh Bolusani\thanks{\texttt{bsuresh@lehigh.edu}}}
\author{Ted K. Ralphs\thanks{\texttt{ted@lehigh.edu}}}
\affil{Department of Industrial and Systems Engineering, Lehigh University,
  Bethlehem, PA, USA}

\maketitle

\begin{abstract}
We describe a framework for reformulating and solving optimization problems
that generalizes the well-known framework originally introduced by Benders. We
discuss details of the application of the procedures to several classes of
optimization problems that fall under the umbrella of multilevel/multistage
mixed integer linear optimization problems. The application of this abstract
framework to this broad class of problems provides new insights and a broader
interpretation of the core ideas, especially as they relate to duality and the
value functions of optimization problems that arise in this context.
\end{abstract}

\graphicspath{{./figures/}}

\section{Introduction \label{sec:intro}}
This paper describes a framework for reformulating and solving optimization
problems that extends the well-known framework of~\cite{bend:parti}. Although
the basic elements of the framework are known, we provide a self-contained
development of the key concepts and illustrate in detail the principles
involved by applying them to the solution of several classes of optimization
problems, including one to which they have not previously been applied. These
classes of problems are all contained under the broad umbrella of what we
informally refer to as \emph{multilevel/multistage mixed integer linear
optimization problems} (MMILPs). MMILPs comprise a broad class of optimization
problems in which multiple decision makers (DMs), with possibly competing
objectives, make decisions in sequence over time. Each DM's decision impacts
the options available to other DMs at other (typically later)
stages\footnote{We use the term ``stage'' in describing the decision epochs of
an MMILP, rather than the alternative ``level'' used in the multilevel
optimization literature because of its broader connotation and connection to
stochastic optimization.}. In economics, these problems fall under the general
umbrella of game theory. We do not formally define the broad class comprising
MMILPs here, but rather describe some specific subclasses contained within it.
Readers wishing to have a more complete overview of MMILPs should refer
to~\cite{BolConRalTah20}.

Although Benders' technique was originally applied to standard mathematical
optimization problems with an underlying structure suggesting a partition of
the variables into exactly \emph{two} sets, it can be similarly applied not
only to more general classes of optimization problems, but by extension, to
$l$-stage problems in which there is an obvious division of the variables into
$l$ sets. Because an $l$-stage MMILP is most naturally defined recursively in
terms of an ($l-1$)-stage MMILP, the very structure of MMILPs seems to suggest
solution by an approach similar to the one suggested by Benders. The recursive
structure also mirrors that of the \emph{polynomial time hierarchy} (PTH,
originally introduced by~\cite{stockmeyer77}), a recursively defined family of
complexity classes into which MMILPs can naturally be categorized. The lowest
level of the PTH is the well-known class $\Pcomplexity$ of problems solvable
in time polynomial in the size of the input, and the $l^{\textrm {th}}$ level
(whose class of primary interest is denoted $\Sigma_l^P$) is comprised of
problems solvable in polynomial time given an oracle for problems in the
$(l-1)^{\textrm {st}}$ level.
The decision versions of MMILPs with $l$ levels are prototypical complete
problems for $\Sigma_l^P$~\citep{jeroslow85}, meaning that all other problems
in the class can be reduced to MMILPs.

Benders' framework is, first and foremost, a technique for reformulation.
Using this technique, MMILPs can be recast as standard mathematical
optimization problems. The reformulation usually results in an exponential
increase in size relative to the original formulation, and a number of
additional transformations may be necessary to get the final problem into a
form in which a blackbox solver can digest it.
Because of the exponential increase in size, the reformulation must generally
be solved either by an approach based on a convergent iterative approximation
scheme or by utilizing a relaxation to obtain bounds that can then be used to
drive a branch-and-bound algorithm. These two approaches are closely related,
as the relaxations required in the latter approach can be obtained by
terminating the iterative approximation procedure before convergence. We
further discuss approaches based on branch and bound in
Section~\ref{subsec:existing-literature-relationships}.

In the remainder of the paper, we focus only on the iterative approximation
approach, which can be seen as a generalization of the cutting-plane method for
solving mixed integer linear optimization problems (MILPs). This approach
can be applied recursively, essentially decomposing the problem by stage, with
the subproblem (introduced formally in Section~\ref{sec:gbPrinciple}) that arises
when solving an $l$-stage problem being an (lexicographic)
optimization problem with $l-1$ stages. The main contributions of this work
are (1) the development of an abstract framework for generalizing the
principles of Benders' technique for reformulation that encompasses
non-traditional problem classes, (2) the specification of an associated
algorithmic procedure that generalizes the standard cutting-plane algorithm
and is based on iterative approximation of functions arising from the
projection of the original problem into the space of a specified subset of
variables, and (3) its application to the solution of \emph{mixed integer
two-stage/bilevel linear optimization problems} (MIBLPs). To our knowledge,
this is the first algorithm for MIBLPs that utilizes a generalized Benders'
approach.
This paper does not aim at discussing efficiency or comparing the algorithms
described herein to alternatives. While we have implemented a proof-of-concept
for the algorithms described here, a full-featured, efficient implementation
would require substantial additional development.

The paper is organized as follows. In
Section~\ref{sec:gbPrinciple}, we discuss the principles underlying our
generalized Benders' framework at a high level in the context of a general
optimization problem, including concepts of bounding functions and general
duality. We also highlight the relationships between this framework and
certain existing algorithms at an abstract level. In
Section~\ref{sec:existingApplications}, we illustrate these principles
concretely with two examples, summarizing existing algorithms for the case in
which there are two stages and the objective functions are the same in both
stages. Section~\ref{sec:MIBLPs} goes into more detail in describing an
algorithm for general MIBLPs, the special case of MMILPs in which there are
only two stages. Finally, we conclude in Section~\ref{sec:conclusions} by
briefly discussing a conceptual extension of the algorithm for MIBLPs to general
MMILPs with $l$ stages.

\section{Benders' Principle \label{sec:gbPrinciple}}

In this section, we introduce the basic principles of the framework. We first
describe it in a very general context and then focus on the special case in
which the objective and constraint functions are \emph{additively separable}.
The idea of such a generalization of Benders' original algorithm is not new.
As far back as the 1970s, \cite{geof:gene} had already proposed a similar
idea. Its application to MMILPs, however, provides new insights and broader
interpretations of the core ideas.

What we generally mean by a Benders-type approach is a technique for
reformulation and/or solution of an optimization problem that operates in a
subspace associated with a specified subset of variables from the original
compact formulation. We refer to this subset of variables as
the \emph{first-stage} variables throughout the paper, although Benders'
technique only really provides a separation of the problem into independent
``stages'' once we assume additive separability. The essential element
underlying any Benders-type method is a \emph{projection} operation.
Projecting an optimization problem means projecting both its feasible region
and its objective function, in the fashion we describe, in order to obtain a
valid reformulation involving only first-stage variables. By ``valid
reformulation,'' we mean one in which the set of optimal solutions of the
projected problem is the projection of the set of optimal solutions of the
original problem, though one could define the concept of validity in other
ways.

The projection operation is natural in applications where the optimal values
of the first-stage variables are of primary concern, while the remaining
variables are present only to model the later-stage effects of the first-stage
decisions. The goal of the projection operation, however, is purely
pragmatic---it is to construct a reformulation that is somehow algorithmically
advantageous. The advantage may either be because the associated relaxations
are more effective or simply because the reformulation has a form that makes
the employment of existing blackbox solvers easier. The reformulation process
necessarily introduces complex functions of the first-stage variables, which
model the effects mentioned above. Algorithms for solving these reformulations
generally construct approximations of these functions, as we detail in the
following sections.

\subsection{General Optimization Problems \label{sec:GP}}

We first consider the following very general form of optimization problem in
which the variables are partitioned into two sets, the \emph{first-}
and \emph{second-stage variables}, denoted by $x \in \Re^{n_1}$ and
$y \in \Re^{n_2}$, respectively. The problem is then
\begin{equation}\label{eq:GP}
	\min\limits_{x \in X,\ y \in Y} \left\{f(x, y) \midd F(x,
  y) \geq 0\right\},
\end{equation}
where $f: \Re^{n_1} \times \Re^{n_2} \rightarrow \Re$ is the objective
function and $F: \Re^{n_1} \times \Re^{n_2} \rightarrow \Re^m$ is the
constraint function, with $X \subseteq \Re_+^{n_1}$ and $Y \subseteq
\Re_+^{n_2}$ denoting the additional disjunctive constraints on the values of
the variables. Typically, we have $X = \mathbb{Z}^{r_1}_+ \times
\mathbb{R}^{n_1-r_1}_+$ and $Y = \mathbb{Z}^{r_2}_+ \times
\mathbb{R}^{n_2-r_2}_+$, and we therefore consider that form of sets for the
remainder of the paper. By convention, we take the optimal objective value to be
$\infty$ if the feasible region
\begin{equation*}
  \Fmcal = \left\{(x, y) \in X \times Y
  \;\middle|\; F(x, y) \geq 0\right\}
\end{equation*}
is empty and $- \infty$ if the problem~\eqref{eq:GP} is unbounded. We assume
that in all other cases, the problem has a finite minimum that can be
attained.

\subsubsection{Projection and the Subproblem}

The simple yet fundamental idea is that~\eqref{eq:GP} can be
equivalently formulated
as
\begin{equation*}
	\min\limits_{x \in X} \left\{ \min\limits_{y \in Y} \left\{f(x, y)
  \midd F(x, y) \geq 0\right\}\right\}.
\end{equation*}

By replacing the inner optimization problem with a function, we obtain the
reformulation
\begin{equation}\label{eq:GP-Proj}
	\min\limits_{x \in X} \phiSS(x),
\end{equation}
in terms of only the first-stage variables, where
\begin{equation}\label{eq:GP-SP}
	\phiSS(x) = \min\limits_{y \in Y} \left\{f(x, y) \midd F(x, y)
  \geq 0\right\} \quad \forall x \in \Re^{n_1}.
\end{equation}

In this new formulation, $\phiSS$ is a function that returns the objective
function value of the optimal feasible combination of values for both first- and
second-stage variables, given fixed values for the first-stage variables.

Although the formulation~\eqref{eq:GP-Proj} does not explicitly involve
projection, we define by convention that $\phiSS(x) = \infty$ if
$x \not\in \proj_x(\Fmcal)$, where
\begin{equation}\label{eq:GP-FR-Proj}
\proj_x(\Fmcal) = \{x \in X \mid F(x, y) \geq 0 \text{ for some } y \in Y\}
\end{equation}
is the projection of the feasible region of~\eqref{eq:GP} into the space of
the first-stage variables. (We similarly define $\phiSS(x) = -\infty$ if the
optimization problem on the right-hand side of~\eqref{eq:GP-SP} is unbounded.)
This means that $\phiSS$ plays a dual role. First, it effectively prevents any
first-stage solution that is not in the projected feasible region from being
considered (provided the projected feasible region is non-empty). Second, it
is also what we earlier described informally as the projection of the original
objective function, since it ensures that the objective function value in the
projected optimization problem with respect to $\hat{x} \in \proj_x(\Fmcal)$
is exactly the value that would have been obtained if solving the original
problem with the first-stage variables fixed to $\hat{x}$. Overall, the
reformulation \eqref{eq:GP-Proj} can then be considered to be the projection
of~\eqref{eq:GP} into the space of the first-stage variables. The evaluation
of $\phiSS$ for particular first-stage solutions is the
aforementioned \emph{subproblem}. More details about its role in the overall
solution process are provided in Section~\ref{sec:overall}.

\subsubsection{Bounding Functions and the Master Problem}

In principle, the optimal value of~\eqref{eq:GP}, as well as an optimal
first-stage solution, can be obtained by solving~\eqref{eq:GP-Proj}. However,
we usually do not have a closed-form description of $\phiSS$ and even when
such closed form exists in theory, its description is typically of
exponential size and would thus be impractical. We therefore replace
$\phiSS$ in~\eqref{eq:GP-Proj} with a dual bounding function
$\uphiSS$ (defined below) to obtain the relaxation known as the \emph{master problem}.
\begin{definition}[\textbf{Dual Bounding Function}]\label{def:LBF}
A function $\uphiSS : \Re^{n_1} \rightarrow \Re \cup \{\pm \infty\}$ is said
to be a dual bounding function with respect to the projection of the
objective function $f$ if
\begin{equation}
\uphiSS(x) \leq f_x(x) \quad \forall x \in \Re^{n_1}.\nonumber
\end{equation}
It is called \emph{strong} at $\hat x \in X$ if
\begin{equation}
\uphiSS(\hat x) = f_x(\hat x). \nonumber
\end{equation}
\end{definition}

Given a dual bounding function, the master problem is then 
\begin{equation}\label{eq:GP-MP-LB}
	\min\limits_{x \in X} \uphiSS(x).
\end{equation}

Naturally, for any relaxation-based method to be practical, solving the
relaxation (in this case, \eqref{eq:GP-MP-LB}) should be easier than solving
the original problem (in this case, \eqref{eq:GP-Proj}). The difficulty of
solving~\eqref{eq:GP-MP-LB}, however, is directly related to the structure of
the function $\uphiSS$ itself. In the cases discussed later, this
function is piecewise linear and the master problem can be formulated as an
MILP.

\subsubsection{Overall Algorithm \label{sec:overall}}

The overall method is to iteratively improve the master problem formulation by
strengthening $\uphiSS$. In iteration $k$, candidate solution
$x^k$ is generated by solving the current master problem (yielding a lower
bound on the optimal value of~\eqref{eq:GP}), and $\phiSS(x^k)$ is then
evaluated (yielding an upper bound on the optimal value of~\eqref{eq:GP}).
The algorithm alternates between solution of
the master problem and the subproblem until upper and lower bounds are equal.

Although the evaluation of $\phiSS(x^k)$ apparently involves only the
determination of the optimal solution \emph{value}, the solution of the
subproblem also typically produces, as a byproduct, a primal-dual proof of
optimality for the problem on the right-hand side of~\eqref{eq:GP-SP}. It is from
this primal-dual proof that we extract a dual bounding function $\uphiSS^k$ that
is strong at $x^k$. The form of this primal-dual proof and the structure of
$\uphiSS^k$ depends strongly on the form of the original problem. We examine
particular cases in Section~\ref{sec:duality}.

It is possible that an algorithm following this general outline will either
converge to a local optimum or not converge at all (see \cite{sahi:conv}),
but the convergence of the method to a global optimal solution can be
guaranteed under two conditions that are
satisfied in many important cases. The first of these is that we update
$\uphiSS$ in each iteration $k$ in such a way that we guarantee
that it is strong not only at $x^k$ but also at all
$x^i$, $i < k$. This can be most easily accomplished by taking the maximum
of the bounding functions generated at each iteration. That is, after
iteration $k$,
\begin{equation} \label{eq:LBF}
  \uphiSS(x) = \max_{1 \leq i \leq k} \uphiSS^i(x),
\end{equation}
where $\uphiSS^i$ is the dual bounding function obtained in
iteration $i \leq k$ of the algorithm. In such cases, the master
problem~\eqref{eq:GP-MP-LB} is usually reformulated using a standard trick to
eliminate the maximum operator by introducing an auxiliary variable $z$ to
obtain 
\begin{equation}\label{eq:GP-MP}
\begin{aligned}
	&& \min\limits_{x \in X} &\; z\\
        & & \text{s.t.} & \; z \geq \uphiSS^i(x) \quad 1 \leq i \leq k.
\end{aligned}
\end{equation}

The formulations~\eqref{eq:GP-MP} and~\eqref{eq:GP-MP-LB} are equivalent in
this case because $z$ must be equal to the maximum of the individual bounding
functions at optimality. In the literature, the constraints
$z \geq \uphiSS^i(x)$ for $1 \leq i \leq k$ in~\eqref{eq:GP-MP}
are sometimes called \emph{Benders' optimality constraints}.
Depending on how the master problem is reformulated, it may also sometimes be
necessary to explicitly exclude $x^k \not\in \proj_x(\Fmcal)$ from the
feasible region of the master problem, in which case the associated constraints
are called \emph{Benders' feasibility constraints}.

The overall approach is outlined in Figure~\ref{fig:gpBenders}.
\begin{figure}
{\bf{Generalized Benders' Framework for Solving}~\eqref{eq:GP}}
\\
\hrule\smallskip
\textbf{Step 0. Initialize} $k \leftarrow 1$, $\textrm{UB}^0 = \infty$, $\textrm{LB}^0 =
-\infty$, $\uphiSS^0(x) = -\infty$ for all $x \in
\Re^{n_1}$. \smallskip\\
\textbf{Step 1. Solve the master problem (lower bound)}
\begin{itemize}
\item[a)] Construct the dual bounding function $\uphiSS(x) =
  \max\limits_{0 \leq i \leq k-1} \uphiSS^i(x)$ and formulate the master
  problem~\eqref{eq:GP-MP}. \smallskip
\item[b)] Solve~\eqref{eq:GP-MP} to obtain an optimal solution
  $(x^{k},\ z^{k})$. Set $\textrm{LB}^k \leftarrow z^{k}$.
\end{itemize}
\textbf{Step 2. Solve the subproblem (upper bound)}
\begin{itemize}
\item[a)] Solve~\eqref{eq:GP-SP} for the given $x^{k}$ to obtain an
  optimal solution $y^{k}$ and strong dual bounding function
  $\uphiSS^{k}$ such that $\uphiSS^k(x^{k}) =
  \phiSS(x^{k})$. Set $\textrm{UB}^k \leftarrow \phiSS(x^{k})$. \smallskip
\item[b)] Termination check: $\textrm{UB}^k = \textrm{LB}^k$? \smallskip
\begin{enumerate}
\item If yes, STOP. $(x^{k}, y^{k})$ is an optimal solution
  to~\eqref{eq:GP}. \smallskip
\item If no, set $k \leftarrow k+1$ and go to Step 1.
\end{enumerate}
\end{itemize}
\hrule\medskip
\caption{Outline of the generalized Benders' decomposition
  framework \label{fig:gpBenders}} 
\end{figure}
Theorem~\ref{thm:genBendersConvergence} due to \cite{hook:logic} shows that
Algorithm~\ref{fig:gpBenders} converges in a finite number of steps under the
additional condition that $\proj_x(\Fmcal)$ is finite.
\begin{theorem}[\cite{hook:logic}, Theorem 2]\label{thm:genBendersConvergence}
If the function $\uphiSS$ is defined as in~\eqref{eq:LBF}
and $\uphiSS^i$ is strong at $x^i$ in each iteration $i$, then
$\uphiSS$ remains a valid dual bounding function that is strong
at all previous iterates and the method converges to the optimal value in a
finite number of iterations under the assumption that
$|\proj_x(\Fmcal)| < \infty$.
\end{theorem}

The proof of this result is rather straightforward. In fact, a slightly more
general result also holds, since the overall dual bounding function need not
be constructed in this particular way, as long as we can ensure that it is
strong at all the previous iterates. In practical implementations, however,
taking the maximum of previous bounding functions is a natural approach and it
is the one we adopt here.

This general framework still leaves many steps unspecified and raises many
questions regarding implementation in specific cases. These questions will be
answered in detail for the several cases of our interest in
Sections~\ref{sec:existingApplications} and~\ref{sec:MIBLPs}.

\subsection{Additively Separable Optimization Problems \label{subsec:GPwSF}}

We now move to the more specific setting that is central to the application of
Benders' method to optimization problems in which the constraint and objective
functions are additively separable.

\begin{definition}[\textbf{Additively Separable Function}]
A function $f: \Re^{n_1} \times \Re^{n_2} \rightarrow \Re$ is \emph{additively
separable} if $\exists\; g: \Re^{n_1} \rightarrow \Re$ and $h: \Re^{n_2}
\rightarrow \Re$ such that $f(x, y) = g(x) + h(y)$ for all $(x, y) \in
\Re^{n_1} \times \Re^{n_2}$.
\end{definition}

When the functions $f$ and $F$ are additively separable, such
separability allows us to reformulate these problems in ways that enhance
intuition and also ease implementation. As such, let $g: \Re^{n_1} \rightarrow
\Re$, $h: \Re^{n_2} \rightarrow \Re$, $G: \Re^{n_1} \rightarrow \Re^m$, and $H:
\Re^{n_2} \rightarrow \Re^m$
be such that $f(x, y) = g(x) + h(y)$ and $F(x, y) = G(x) + H(y)$ for all $(x,
y) \in \Re^{n_1} \times \Re^{n_2}$. Because we are specifically interested in
the case of linear functions, we also introduce a right-hand side $b
\in \Re^m$, as is standard for problems involving linear functions. We
then obtain the new form of general optimization problem
\begin{equation}\label{eq:GP-AS}
	\min\limits_{x \in X,\ y \in Y} \left\{g(x) + h(y) \midd G(x) +
  H(y) \geq b\right\}
\end{equation}
that we consider in the rest of the paper.

\subsubsection{Projection and the Value Function}

A reformulation of~\eqref{eq:GP-AS} analogous to~\eqref{eq:GP-Proj}, obtained
upon projecting into the space of the first-stage variables, is
\begin{equation}\label{eq:GP-VF}
	\min\limits_{x \in X} \left\{g(x) + \phiVF\left(b -
  G(x)\right)\right\},
\end{equation}
where 
\begin{equation}\label{eq:SP-VF}
	\phiVF(\beta) = \min\limits_{y \in \Re_+^{n_2}} \left\{h(y)
  \midd y \in \Pcal_2(\beta) \cap Y \right\} \quad \forall \beta \in \Re^{m},
\end{equation}
and
\begin{equation*}
\P_2(\beta) = \left\{y \in \Re_+^{n_2} \midd H(y) \geq \beta\right\} \quad
\forall \beta \in \Re^m 
\end{equation*}
is a parametric family of polyhedra containing points satisfying the
second-stage feasibility conditions, which can now be considered fully
independently, due to the additive separability. By convention, $\phiVF(\beta)
= +\infty$ for $\beta \in \Re^{m}$ if $\Pcal_2(\beta) \cap Y = \emptyset$, and
$\phiVF(\beta) = -\infty$ for $\beta \in \Re^{m}$ if the problem on the
right-hand side of~\eqref{eq:SP-VF} is unbounded.

As opposed to the earlier-defined function $\phiSS$, which was
parameterized on the first-stage solution, $\phiVF$ is parameterized on the
right-hand side of the associated second-stage optimization problem (which is
in turn determined by the first-stage solution). The second-stage optimization
problem, analogous to the earlier defined subproblem~\eqref{eq:GP-SP}, is to
evaluate $\phiVF$ at a specific right-hand side.
In the context of the framework described in Figure~\ref{fig:gpBenders},
$\phiVF$ would be evaluated in iteration $k$ of the algorithm at the
right-hand side $\beta^k = b - G(x^k)$.

Those readers familiar with the more general duality theory associated with
mixed integer linear optimization problems (see,
e.g., \cite{nemhauser1988integer} and \cite{GuzRal07}) will recognize $\phiVF$ as the
\emph{value function} of the second-stage problem. The value function of the
second-stage problem and the associated dual problem are crucial elements of
the framework in the additively separable case, so we now briefly review these
basic concepts.

\subsubsection{General Duality and Dual Functions \label{sec:duality}}

The function generally referred to as ``the'' value function of an
optimization problem is one that returns the optimal objective value for a given
right-hand side vector. As in the general framework from Section~\ref{sec:GP},
we form the master problem by replacing this value function with a function
that bounds it from below. Such functions are known as
\emph{dual functions}, so called because they can be interpreted as solutions
to a general \emph{dual problem} and reflect the essential role of duality in
Benders-type reformulations. For a particular right-hand side $\hat\beta \in
\Re^m$, the general dual problem
\begin{equation}\label{eq:SP-VF-GD}
\max\limits_{D \in \Upsilon^m}
\left\{
D(\hat\beta) \midd D(\beta) \leq \phiVF(\beta) \quad \forall \beta \in
\Re^{m} 
\right\}
\end{equation}
associated with~\eqref{eq:SP-VF} is an optimization problem over a class
$\Upsilon^m \subseteq \left\{\upsilon \mid \upsilon : \Re^m \rightarrow
\Re\right\}$ of real-valued functions. The objective of the dual problem is
to construct a function that bounds the value function from below and for
which the bound is as strong as possible at $\hat\beta$. As such, we define a
\emph{(strong) dual function} as follows.
\begin{definition}[\textbf{Dual Function}]\label{def:dualfn}
A dual function $D : \Re^{m} \rightarrow \Re \cup \{\pm \infty\}$ is one that
satisfies $D(\beta) \leq \phiVF(\beta)$ for all $\beta \in \Re^{m}$.  It is
called \emph{strong} at $\hat \beta \in \Re^m$ if $D(\hat \beta) =
\phiVF(\hat \beta)$.
\end{definition}

The dual problem itself is called \emph{strong} if $\Upsilon^m$ is guaranteed
to contain a strong dual function. As long as the value function itself is
real-valued\footnote{When the value function is not real-valued everywhere, we
have to show that there exists a real-valued function that coincides with the
value function whenever the value function is real-valued and is itself
real-valued everywhere else, but is still a feasible dual function
(see~\cite{wolsey81}).} and is a member of $\Upsilon^m$, then the dual problem
will be strong, since the value function itself is an optimal solution
of~\eqref{eq:SP-VF-GD}.

Exact solution algorithms that produce certificates of optimality typically do
it by providing a primal solution, which certifies an upper bound, and a dual
function (solution to~\eqref{eq:SP-VF-GD}), which certifies a lower bound.
When these bounds are equal, the combination provides the required certificate
of optimality. Dual functions can be obtained in a variety of ways, but one
obvious way to construct them is to consider the value function of a
relaxation of the problem. Most solution algorithms for linear and mixed
integer linear optimization problems work by iteratively constructing such a
dual function.

The connection between the general dual and Benders' framework should be
clear. The strong dual function constructed as a certificate of optimality
when solving the subproblem (evaluating $\phiVF$) is a function that can
be directly used in strengthening the global dual function that defines the
current master problem. In fact, it is useful to think of the subproblem not
as that of evaluating $\phiVF$, but rather of solving a dual problem of the
form~\eqref{eq:SP-VF-GD} to obtain a strong dual function, which is what we
actually need for forming the master problem in the next iteration.

\subsubsection{Overall Algorithm}

The overall method is largely similar to that described in
Figure~\ref{fig:gpBenders}. We relax the reformulation~\eqref{eq:GP-VF} to
obtain a master problem
\begin{equation}\label{eq:MP-VF}
\begin{aligned}
	&& \min\limits_{x \in X} &\; g(x) + z\\
        & & \text{s.t.} & \; z \geq \uphiVF\left(b - G(x)\right),
\end{aligned}
\end{equation}
in which the value function $\phiVF$ is replaced by a dual function $\uphiVF$.
Solving this master problem in iteration $k$ yields a solution $(x^k,
z^k)$ and a lower bound $g(x^k) + z^k$. We then evaluate $\phiVF$ at $b -
G(x^k)$ to obtain a dual function $\uphiVF^k$ strong at $b -
G(x^k)$ and an upper bound $g(x^k) + \phiVF(b - G(x^k))$. This dual function
is combined with previously produced such functions to obtain a global dual
function that is strong at all previous iterates, ensuring eventual
convergence under the same conditions as in
Theorem~\ref{thm:genBendersConvergence}.

Naturally, the exact form and structure of the dual functions involved is
crucially important to the tractability of the overall algorithm, as we do
need a method of (re)formulating and solving the master problem in each
iteration. In the cases discussed in this paper, the dual function takes on
relatively simple forms. For linear optimization problems (LPs), the value function
is convex and there is always a strong dual function that is a simple linear
function. This linear function is an optimal solution to the usual LP dual,
which is a subgradient of the LP value function.

In the case of mixed integer linear optimization, dual functions can be
obtained as a by-product of a branch-and-bound algorithm. Roughly speaking,
the lower bound produced by a branch-and-bound algorithm is the minimum
of lower bounds produced for
the individual subproblems associated with the leaf nodes of the
branch-and-bound tree. Thus, the overall dual function is the minimum of dual
functions for these subproblems. In the MILP case, the subproblem dual
functions utilized are affine functions derived from the dual of the LP
relaxation of the subproblem associated with a given node. Thus, in the
simplest case, the dual function is the minimum of affine functions.

This method of constructing dual functions from branch-and-bound trees can be
extended to virtually any problem that can be solved by a relaxation-based
branch-and-bound algorithm. The lower bound arising from the branch-and-bound
tree is the minimum of lower bounds on individual subproblems, which are
typically (but not always) derived as dual functions of convex relaxations.
The overall dual function is thus a minimum of dual functions for individual
leaf nodes, as in the MILP case. In Sections~\ref{subsubsec:2SSMILPs-proj}
and~\ref{sec:reaction-fcn} below, we describe in detail the application of
this principle to the derivation of dual functions for the MILP and
lexicographic MILP cases.

\subsection{Relationship to Other Methodologies
\label{subsec:existing-literature-relationships}}

An important question is how general this framework is and how
existing algorithms are related. The framework presented is at a high level of
abstraction and, therefore, very generic. It is likely that almost all
Benders-type procedures can be seen as special cases or at least subtle
variations on the same theme, since any reformulation of the problem in a
subspace would necessarily involve a projection operation in some form. Of
course, as with all algorithms described at this level of abstraction, many
important details need to be filled in to attain a practical algorithm for
specific problem classes, so the development of customized algorithms for
specific classes is still very much needed. The framework only provides a way
to understand relationships and a starting template for building algorithms.

As mentioned earlier, Benders' reformulation technique generally increases the
size of the formulation exponentially relative to the original compact
formulation and the reformulation may well be in a wholly different class than
the original problem (e.g., the projection of a linear problem may involve
nonlinear functions, see Section~\ref{subsec:LPIllustration}). This is
analogous to the way in which a minimal description of the convex hull of
feasible solutions of a standard MILP has a description of exponential size
with respect to its original formulation.
With a reformulation of exponential size, as in the solution of MILPs, the most
obvious approaches are ones based on iterative outer approximation, suggesting
a generalized cutting plane-type method. Such methods produce a sequence of
relaxations similar to what we referred to as the ``master problem,'' whose
solutions converge to the exact optimum. From this point of view, a
traditional cutting-plane algorithm for solving the problem in its original
compact form is just one possible approach on a continuum ranging from
projecting out none of the variables to projecting out all the variables. (The
latter option of projecting out all the variables would essentially be
equivalent to constructing the value function of the original problem.)

To take a specific example,~\cite{codatoFischetti06} propose a Benders-type
algorithm that employs so-called \emph{combinatorial Benders' cuts} (also
known as no-good cuts). These cuts can be seen as Benders' feasibility
constraints in the framework of this paper because they remove individual
first-stage solutions that are not in the projection of the feasible region.
The set of all such cuts, along with the requirement that the first-stage
variables be binary, provides an exact description of the projected feasible
region. And since the objective function does not depend on the second-stage
variables, Benders' optimality constraints are not needed. As another example,
\cite{sen2006decomposition} discuss a Benders-type algorithm for solving
two-stage stochastic mixed integer linear optimization problems (2SSMILPs).
They propose linear Benders' cuts that can be viewed as Benders' optimality
constraints. These cuts are obtained by applying the
disjunctive cut principle \citep{balas:79} to the disjunction based on LP
relaxations of leaf nodes of the branch-and-bound tree for solving MILP
subproblems.

As with cutting plane methods for MILPs, when convergence of such a method is
slow, the process can be terminated early to obtain a bound. This bounding
procedure can be embedded within a branch-and-bound framework to obtain what
is called a \emph{branch-and-Benders'-cut} algorithm in the literature and is
essentially a generalized branch-and-cut method. \cite{RahCraGenRei17} provides
a list of works that propose such an algorithm.
Note that another reason for embedding the procedure in a branch-and-bound
framework is to simplify the master problem by restricting the form of the
functions involved so that exact reformulation is no longer possible. For
example, we may want the master problem to be convex for reasons of
efficiency. Because the functions involved can be non-convex even when the
original formulation contains only linear functions, this would rule out exact
reformulation. By partitioning the feasible region into small enough regions,
convex reformulations become possible. This is analogous to the notion of
spatial branching in mixed integer nonlinear optimization.

Hybrid options are also possible, especially in multistage optimization. The
most obvious approach is to allow the later-stage variables to remain in the
master problem. This obviates the need for dual bounding functions and
avoids some of the difficulties associated with projection of the feasible
region, while retaining the possibility of viewing the optimization as being
over the space of the first-stage variables only. This is the approach taken
in most branch-and-cut algorithms for solving bilevel optimization
problems~\citep{TahRalDeN20}. We discuss this possibility further in the
context of MIBLPs in the literature review of Section~\ref{sec:MIBLPs}.

\section{Applications From the Literature \label{sec:existingApplications}}

In this section, we describe the application of the generalized framework
presented in Section~\ref{sec:gbPrinciple} to derive algorithms already
existing in the literature. We describe these applications here to emphasize
their commonality, and to provide concrete examples in settings in which the
application of the principles is relatively straightforward and the
abstractions reduce to familiar algorithmic concepts.

\subsection{Linear Optimization Problems \label{subsec:LPIllustration}}

We begin by considering the application of Benders' framework to the standard LP
\begin{equation}\label{eq:LP}
\min \left\{c^\top x + d^\top y \midd Ax + Gy \geq b,\; (x, y) \in \Re_+^{n_1} \times
\Re_+^{n_2} \right\},
\end{equation}
where $A \in \mathbb{Q}^{m \times n_1}$, $G \in \mathbb{Q}^{m \times n_2}$, $b
\in \mathbb{Q}^m$, $c \in \Q^{n_1}$, and $d \in \Q^{n_2}$. This problem is the
special case of~\eqref{eq:GP-AS} in which $g(x) = c^\top x$, $h(y) = d^\top y$, $G(x) =
Ax$, and $H(y) = Gy$ for all $x \in X = \Re^{n_1}_+$ and $y \in Y =
\Re^{n_2}_+$.

\subsubsection{Projection}
Projecting into the space of
the first-stage variables, we obtain the reformulation
\begin{equation}\label{eq:LP-Proj-VF}
\min\limits_{x \in \Re_+^{n_1}} \left\{c^\top x + \phi_{LP}(b - Ax)
\right\},
\end{equation}
where 
\begin{equation*}
\phi_{LP}(\beta) = \min \left\{d^\top y \midd Gy \geq \beta,\; y \in \Re_+^{n_2}
\right\} \quad \forall \beta \in \Re^m,
\end{equation*}
is the value function of the second-stage optimization problem, which is an LP.
This reformulation is nothing more than the instantiation of the
reformulation~\eqref{eq:GP-VF} in the context of~\eqref{eq:LP}.

Value functions are well-studied and well-understood in the linear
optimization case (see, e.g., \cite{bertsimas1997introduction} for details).
The structure of $\phi_{LP}$ arises from that of the feasible region 
\begin{equation*}
\Dmcal = \left\{\eta \in \Re^m_+ \midd G^\top \eta \leq d\right\}
\end{equation*}
of the standard LP dual of the second-stage problem, which is the LP
\begin{equation} \label{eq:LP-Dual}
  \max_{\eta \in \Dmcal}\ \eta^\top (b - A\hat{x})
\end{equation}
when the first-stage solution is $\hat{x} \in \Re^{n_1}_+$. Although it may
not be obvious, this dual problem is precisely equivalent to the general
dual~\eqref{eq:SP-VF-GD} in the LP case. This can be seen by noting the
constraints ensure that the dual solution is a subgradient of $\phi_{LP}$ and
hence represents a (linear) dual function. By noting that the above maximum
can be taken over the set $\Emcal$ of extreme points of $\Dmcal$, assuming
$\Dmcal$ is bounded, it is easy to derive that
\begin{equation}\label{eq:LP-VF-Struct}
\phi_{LP}(\beta) = \max\limits_{\eta \in \mathcal{E}} \left\{\eta^\top \beta \right\} \quad
\forall \beta \in \Re^m.
\end{equation}

That is, the value function is the maximum of linear functions corresponding
to members of $\Emcal$. Although this function is convex and nicely
structured, the cardinality of $\Emcal$ is exponential in general, so
enumerating them is impractical. The global dual function we use to construct
the master problem is thus formed from a small collection of these extreme
points, as described next.

\subsubsection{Master Problem}
In accordance with the principles described
earlier, we form the master problem by replacing $\phi_{LP}$
in~\eqref{eq:LP-Proj-VF} with a global dual function $\uphi_{LP}$
that is the maximum of the strong dual functions produced in each iteration of
the algorithm. In this context, the strong dual functions produced at each
iteration are the linear functions associated with solutions to the
dual~\eqref{eq:LP-Dual} of the second-stage problem. This results in the
master problem
\begin{equation}\label{eq:LP-MP}
\begin{aligned}
	& & \min & \; c^\top x + z \\
        & & \text{s.t.} & \; z \geq \uphi^i_{LP}(b - A x) = {\eta^i}^\top(b
  - A x), \quad 1 \leq i \leq k\\
	& & & x \in \Re_+^{n_1},\\
\end{aligned}
\end{equation}
in iteration $k$, where $\eta^i \in \Emcal$ is an optimal solution
of~\eqref{eq:LP-Dual} in iteration $i$. Note that if the optimal solution to the
master problem in iteration $k$ is $(x^k, z^k)$, then we have
\begin{equation*}
  z^k = \uphi_{LP}(b - A x^k) = \max_{1 \leq i \leq k}\ {\eta^i}^\top(b - A
  x^k), 
\end{equation*}
as desired, and this master problem is the equivalent to~\eqref{eq:MP-VF} in
this context.

\subsubsection{Overall Algorithm}
In iteration $k$ of the algorithm, we begin by solving a master
problem~\eqref{eq:LP-MP} to obtain its optimal solution $x^k$. The subproblem
is then to evaluate $\phi_{LP}$ at $\beta^k = b - A x^k$ by solving the
dual~\eqref{eq:LP-Dual} to obtain $\eta^k \in \Emcal$. By defining
$\uphi^k_{LP}(\beta) = {\eta^k}^\top \beta$ for all $\beta \in \Re^m$, we
obtain that $\uphi^k_{LP}$ is a dual function for $\phi_{LP}$ that is
strong at $\beta^k$.

The overall method is then to add one constraint of the form
\begin{equation}\label{eq:Opt-Cut}
z \geq {\eta^k}^\top (b - Ax)
\end{equation}
to the master problem in each iteration $k$ in which~\eqref{eq:LP-Dual} has a
finite optimum. In case~\eqref{eq:LP-Dual} in iteration $k$ is unbounded, then
$x^k$ is not a
member of the projection of the feasible region (defined as
in~\eqref{eq:GP-FR-Proj}) and we instead add a constraint of the form
\begin{equation}\label{eq:Feas-Cut}
0 \geq {\sigma^k}^\top (b - A x),
\end{equation}
where $\sigma^k$ is the extreme ray of $\Dmcal$ that proves infeasibility
of the second-stage problem.

Observe that the master problem, the subproblem, and the constraints added
to the master problem are identical to the corresponding components in a
classical Benders' decomposition algorithm for LPs.  Conventionally, in the
context of LPs, these constraints are referred to as Benders' cuts
with~\eqref{eq:Opt-Cut} being \emph{Benders' optimality cut}
and~\eqref{eq:Feas-Cut} being \emph{Benders' feasibility cut}.
Next, we discuss how to further generalize to the class of 2SSMILPs.

\subsection{Two-Stage Stochastic Mixed Integer Linear Optimization
  Problems \label{subsec:2SSMILPs}} 

The case of 2SSMILPs generalizes the case discussed in the previous section in
two important ways. First, we introduce stochasticity, which is modeled by
specifying a finite number of possible scenarios in the second stage, resulting
in a block-structured constraint matrix overall. This generalization on its own
is relatively straightforward and results in the method known in the literature
as the L-shaped method for solving stochastic linear optimization problems with
recourse~\citep{van1969shaped}. However, we also wish to allow integer variables
into the second stage. Although this does not require any modification of the
framework itself, it results in a more complex structure for the value function
of the second-stage problem and hence, a more complex reformulation for the
master problem. We now summarize an algorithm for solving 2SSMILPs that was
originally developed by \cite{HasRal14a}.

To model the stochasticity, we introduce a random variable $U$ over an outcome
space $\Omega$ representing the set of possible scenarios for the second-stage
problem. We assume that $U$ is discrete, i.e., that the outcome space $\Omega$
is finite so that $\omega \in \Omega$ represents which of the finitely many
scenarios is realized. In practice, this assumption is not very restrictive, as
one can exploit any algorithm for the case in which $\Omega$ is assumed finite
to solve cases where $\Omega$ is not (necessarily) finite by utilizing a
technique for discretization, such as \emph{sample average approximation}
(SAA)~\citep{shapiro2003monte}.

Under these assumptions, a 2SSMILP is then a problem of the form
\begin{equation}\label{eq:2SSMILP}
 \begin{aligned}
	\min  &\ c^\top x + {\mathbb E}_{\omega \in \Omega} \;\left[{d^2}^\top y_\omega\right]\\
	\text{s.t.} &\ A^1 x \geq b^1 \\
	&\ A^2_\omega x + G^2y_\omega \geq b^2_\omega \quad \forall
	\omega \in \Omega \\
	&\ x \in X, y_\omega \in Y,
 \end{aligned}
\end{equation}
where $c \in \Q^{n_1}$, $d^2 \in \Q^{n_2}$, $A^1 \in \Q^{m_1 \times n_1}$, $G^2 \in
\Q^{m_2 \times n_2}$, and $b^1 \in \Q^{m_1}$. $A^2_\omega
\in\Q^{m_2\times n_1}$ and $b^2_\omega \in \mathbb{Q}^{m_2}$ represent the
realized values of the random input parameters in scenario $\omega \in \Omega$,
i.e., $U(\omega) = (A^2_\omega, b^2_\omega)$.  The first term in the objective
function represents the immediate cost of implementation of the first-stage
solution, while the second term is an expected cost over the set of possible
future scenarios.

\subsubsection{Projection \label{subsubsec:2SSMILPs-proj}}

We now reformulate the problem by exploiting two important properties. First,
since $U$ is discrete, we may associate with it a discrete probability
distribution defined by $p \in \Re^{|\Omega|}$ such that $0 \leq p_\omega \leq
1$ and $\sum_{\omega \in \Omega} p_\omega = 1$, where $p_\omega$ represents
the probability of the scenario $\omega \in \Omega$. This allows us to replace
the expectation above with a finite sum. Second, we note that the second-stage
problem has a natural block structure so that the problem decomposes
perfectly into $|\Omega|$ smaller problems, which differ only in the right-hand
side vector once the first-stage solution is fixed (an
advantage of using Benders' approach in this setting). Thus, when we project this
problem as before into the space of the first-stage variables, we obtain the
reformulation
\begin{equation} \label{eq:2SSMILP-SE}
\min \left\{c^\top x + \sum_{\omega \in \Omega} p_\omega \phiIP(b^2_\omega -
A^2_\omega x) \midd A^1x \geq b^1,\; x \in X \right\}, 
\end{equation}
where the value function $\phiIP$ associated with the second-stage problem is
defined by
\begin{equation}\label{eq:MILP-VF}
\phiIP(\beta) = \min\limits_{y \in \Re^{n_2}_+} \left\{{d^2}^\top y \midd y \in
\Pcal_2(\beta) \cap Y \right\} \quad 
\forall \beta \in \Re^{m_2},
\end{equation}
and
\begin{equation}\label{eq:MILP-VF-FR}
\P_2(\beta) = \left\{y \in \Re_+^{n_2} \midd G^2 y \geq \beta\right\} \quad
\forall \beta \in \Re^{m_2}.
\end{equation}

By convention, $\phiIP(\beta) = \infty$ if the feasible region $\{y \in
\Re_+^{n_2} \mid y \in \Pcal_2(\beta) \cap Y\}$ is empty, and $\phiIP(\beta) =
-\infty$ if the second-stage problem associated with $\beta$ is unbounded,
which results in $\phiIP(\beta) = -\infty$ for all $\beta \in \Re^{m_2}$. Note
that the single second-stage value function that appeared in the analogous
reformulation in Section~\ref{subsec:GPwSF} is replaced here with expected
value of the value function across all scenarios, expressed as a weighted sum.
Although each scenario results in a separate outcome, the evaluation of those
scenario outcomes is done in principle using the single second-stage value
function $\phiIP$, which links the scenarios. Thus, the subproblem is to
evaluate the weighted sum of $\phiIP$ across all scenarios for a fixed
first-stage solution.

It is clear from~\eqref{eq:2SSMILP-SE} that solving~\eqref{eq:2SSMILP}
(directly or iteratively) requires exploiting the structure of the MILP value
function~\eqref{eq:MILP-VF}, just as we exploited the structure of the LP value
function in solving~\eqref{eq:LP-Proj-VF} in
Section~\ref{subsec:LPIllustration}. Therefore, we now take a quick detour to
discuss this structure and construction of a strong dual in the MILP case.

\paragraph{Value Function. \label{sec:MILPVF}}

The structure of the value function of an MILP is by now well-studied. Early
foundational works include~\cite{johnson73, johnson74, blair1977value,
blair1979value, blair1982value, blair1984constructive, jeroslow78,
jeroslow79_2, bachem80, bank1983non} and \cite{blair1995closed}. A number of
later works
extended these foundational results~\citep{GuzRal07, Guzelsoy2009, HasRal14,
Hassanzadeh2015}. What follows is a summary of results from these later works.

Let us first look at an example to get some intuition about the structure
of $\phiIP$.
\begin{example}
Consider the following parametric MILP, described by its associated value
function, plotted in Figure~\ref{fig:IPVF}.
\begin{equation}\label{eq:ip}
\begin{aligned}
	\phiIP(\beta) & = & \min & \; 2 y_1 + 4 y_2 + 3 y_3 + 4 y_4\\
        && \text{s.t.} & \; 2 y_1 + 5 y_2 + 2 y_3 + 2 y_4 \geq \beta,\\
        && & \; y_1, y_2, y_3 \in \Z_+,\ y_4 \in \Re_+
\end{aligned}
\end{equation}
\begin{figure}
\includegraphics[width=\textwidth]{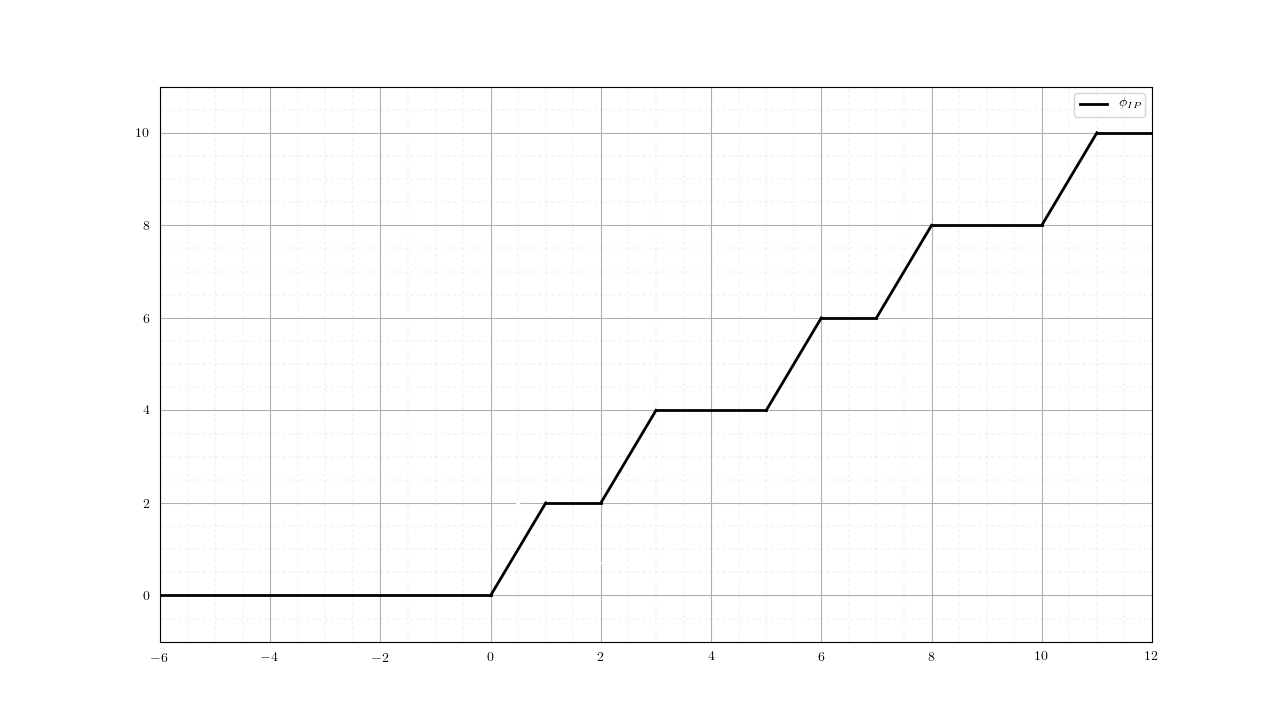}
\caption{MILP value function $\phiIP$~\eqref{eq:ip} \label{fig:IPVF}}
\end{figure}
\qed
\end{example}

The function shown in Figure~\ref{fig:IPVF} is observed to be piecewise
linear, non-decreasing, non-convex, and non-concave. These are all properties
of general value functions of the form~\eqref{eq:MILP-VF},
but there are also several other important properties that are not evident
from this simple example. In particular, the value function may be
discontinuous, but is always lower semi-continuous and also subadditive.

Another important property of $\phiIP$ that is evident from
Figure~\ref{fig:IPVF} is that its epigraph is the union of a set of convex
radial cones. These cones are translations of the epigraph of the value
function of a single parametric LP, the so-called \emph{continuous
  restriction} of the given MILP (defined later in
Section~\ref{sec:reaction-fcn}) resulting from fixing the integer
variables. \cite{HasRal14} further proved that the value function
can be described within any bounded region by specifying a finite set of
points of \emph{strict local convexity} of $\phiIP$, which are the locations of
the extreme points of these radial cones (assuming the epigraph of the LP
value function is a pointed cone). This resulted in a finite discrete
representation of $\phiIP$ (see~\cite{HasRal14} for additional details and
formal results).

Although there exist effective algorithms for evaluating $\phiIP$ for a single
fixed right-hand side $\hat{\beta}$ (e.g., any method for solving the associated
MILP), it is difficult to explicitly construct the entire function because this
evidently requires solution of a sequence of MILPs. Algorithms for evaluating
$\phiIP$ at a right-hand side $\hat{\beta}$, such as branch-and-bound algorithm,
do, however, produce information about its structure beyond the single value
$\phiIP(\hat{\beta})$. This information comes in the form of a dual function that
is strong at $\hat{\beta}$. We next describe the form and structure of these
dual functions.

\paragraph{Dual Functions. \label{sec:MILPVF-dual}}
We focus here on dual functions from the branch-and-bound algorithm, which is
the most widely used solution method for solving MILPs. We refer the reader
to \cite{Guzelsoy2009} and \cite{GuzRal07} for an overview of other methods.
\cite{wolsey81} was the first to propose that dual functions could be
extracted from branch-and-bound trees, as described in the following result.
\begin{theorem}[\cite{wolsey81}, Theorem 19] \label{thm:milp-lbf}
Let $\hat{\beta} \in \Re^{m_2}$ be such that $\phiIP(\hat{\beta}) < \infty$ and
suppose $T$ is the set of indices of leaf nodes of a branch-and-bound tree
resulting from evaluation of $\phiIP (\hat{\beta})$. Then there exists a
dual function $\uphiIP : \Re^{m_2} \rightarrow \Re \cup \{\pm \infty\}$
of the form
\begin{equation} \label{eq:MILP-DF}
\uphiIP (\beta) = \min_{t \in T}\ \left(\beta^\top \eta^t  + \alpha^t\right) \quad
\forall \beta \in \Re^{m_2},
\end{equation}
where $\eta^t \in \Re^{m_2}$ is an optimal solution to the dual of the LP
relaxation associated with node $t$ and $\alpha^t \in \Re$ is the product of
the optimal reduced costs and variable bounds of this LP relaxation. Further,
$\uphiIP$ is strong at $\hat{\beta}$, i.e.,
$\uphiIP(\hat{\beta})= \phiIP(\hat{\beta})$.
\end{theorem}

The interpretation of the function $\uphiIP$ in~\eqref{eq:MILP-DF} is
conceptually straightforward. The solution to the LP relaxation of node $t$ of
the branch-and-bound tree yields the dual function $\beta^\top
\eta^t + \alpha^t$, which bounds the optimal value of the relaxation
problem associated with that node. The overall lower bound yielded by the tree
is then the smallest bound yielded by any of the leaf nodes. This is the usual
lower bound yielded by a branch-and-bound-based MILP solver during the
solution process. Finally, we obtain $\uphiIP$ by interpreting the optimal
solution to the dual of the LP relaxation in each node as a function of
$\beta$. There are additional subtle details involving construction of
appropriate dual functions for infeasible nodes, but we omit these details
here.

In principle, stronger dual functions can be obtained. For example, stronger
functions can be constructed from the branch-and-bound tree by considering
non-leaf nodes, suboptimal dual solutions arising during the solution process,
full LP value function at each leaf node $t$ instead of a single hyperplane
$\beta^\top \eta^t + \alpha^t$, etc. Further details on these methods are mentioned
in \cite{GuzRal07} and \cite{HasRal14a}.

\begin{example}\label{eg:toyMILPLBFs}
Figure~\ref{fig:IPVF2LBFs} shows the dual functions obtained upon applying the
result in Theorem~\ref{thm:milp-lbf} to the MILP~\eqref{eq:ip}. We solve this
MILP with three values of the right-hand side $\beta$.
\begin{figure}
\includegraphics[width=\textwidth]{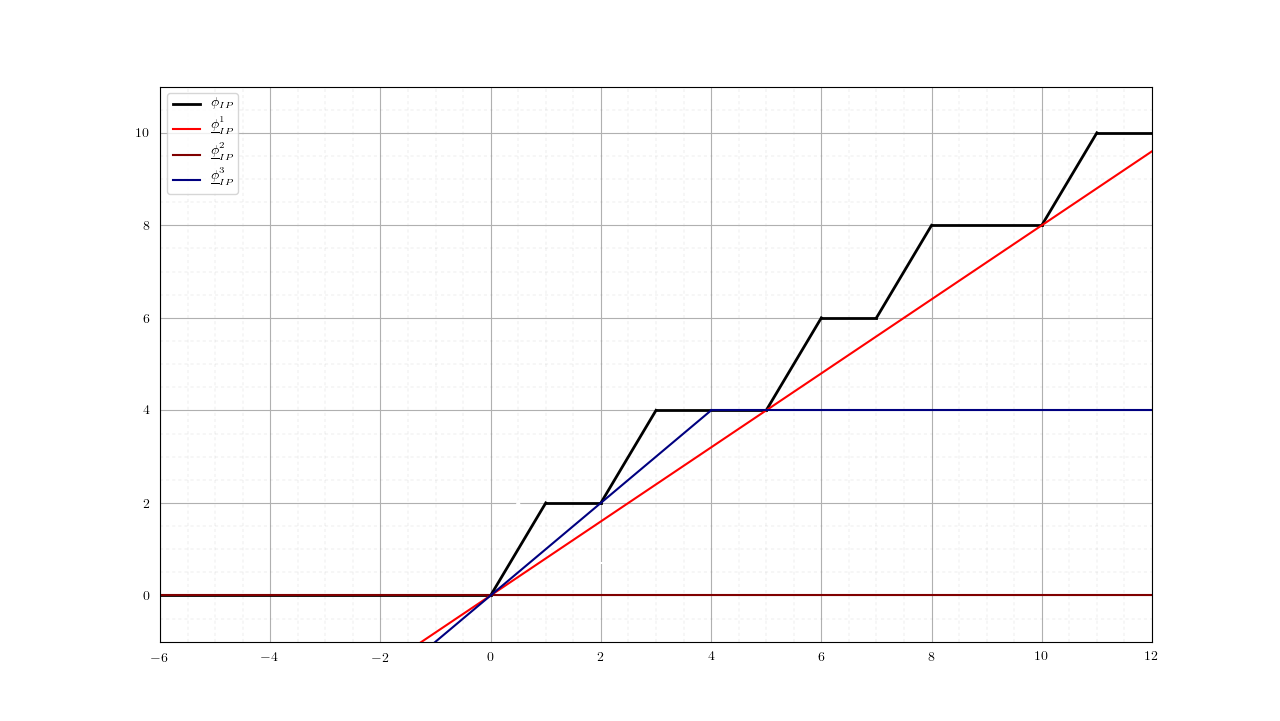}
\caption{Dual functions for~\eqref{eq:ip} \label{fig:IPVF2LBFs}}
\end{figure}
\begin{itemize}
\item $\beta = 5$: There is only one node in the associated branch-and-bound
  tree with the optimal dual solution $\eta = 0.8$, $\underline\eta = (0.4,
    0, 1.4, 2.4)$, and $\overline\eta = (0, 0, 0, 0)$. This results in the dual
    function
\begin{equation*}
\uphiIP^1(\beta) = 0.8 \beta \quad \forall \beta \in \Re.
\end{equation*}
\item $\beta = 0$: There is still only one node in the tree with the optimal
  dual solution $\eta = 0$, $\underline\eta = (2, 4, 3, 4)$, and
    $\overline\eta = (0, 0, 0, 0)$. This results in the dual function
\begin{equation*}
\uphiIP^2(\beta) = 0 \quad \forall \beta \in \Re.
\end{equation*}
\item $\beta = 2$: There are three nodes in the tree, i.e., one root node and
  two leaf nodes resulting from the branching disjunction $y_2 \leq 0 \vee y_2
    \geq 1$. The optimal dual solution and the resulting dual function $\beta^\top
    \eta^t + \alpha^t$ at each leaf node $t$ are mentioned in
    Table~\ref{tab:toyMILPLBFs}.
\begin{table}
\centering
\caption{Data for construction of the dual function from the branch-and-bound
  tree in Example~\ref{eg:toyMILPLBFs} \label{tab:toyMILPLBFs}}
\begin{tabular}{cccccc}
\hline\noalign{\smallskip}
$t$ & Branching constraint & $\eta^t$ & $\underline\eta^t$ & $\overline\eta^t$ &
$\beta^\top\eta^t + \alpha^t$\\
\noalign{\smallskip}\hline\noalign{\smallskip}
1 & $y_2 \leq 0$ & 1 & (0, 0, 1, 2) & (0, -1, 0, 0) & $\beta$ \smallskip\\
2 & $y_2 \geq 1$ & 0 & (2, 4, 3, 4) & (0, 0, 0, 0) & 4 \\
\noalign{\smallskip}\hline
\end{tabular}
\end{table}
This results in the dual function
\begin{equation*}
\uphiIP^3(\beta) = \min \{\beta, 4\} \quad \forall \beta \in \Re.
\end{equation*}
\end{itemize}

Naturally, as in the formulation of the master problem, the
value function approximation can be improved by taking the maximum of multiple
dual functions strong at different right-hand sides. In the above example, the
dual function
$\max\{\uphiIP^1(\beta),\uphiIP^2(\beta), \uphiIP^3(\beta)\}$
  is already seen to be a reasonable approximation of the full value function.
\qed
\end{example}

As mentioned earlier, solving the subproblem in an iteration requires evaluating
$\phiIP$ at $|\Omega|$ right-hand side vectors corresponding to $|\Omega|$
scenarios, for a fixed first-stage solution. Specifically, in iteration $k$ of
the algorithm, we solve
\begin{equation*}
\begin{aligned}
	\phiIP(b^2_\omega - A^2_\omega x^k) = \min & \; {d^2}^\top y \\
        \text{s.t.} & \; G^2y \geq b^2_\omega - A^2_\omega x^k \\
        & \quad y \in Y
\end{aligned}
\end{equation*}
for all $\omega \in \Omega$, where $x^k$ is the fixed first-stage solution
in the current iteration. The result is a scenario dual function
${\uphiIP^k}_{\omega}$ of the form~\eqref{eq:MILP-DF} for each $\omega \in
\Omega$.

\subsubsection{Master Problem}

By exploiting the specific structure of the dual functions described in the
previous section, we can straightforwardly adapt the algorithmic framework from
Section~\ref{sec:gbPrinciple} to obtain an algorithm for
solving~\eqref{eq:2SSMILP-SE} similar to that derived by the authors
in~\cite{HasRal14a}.

Introducing auxiliary variables $z_\omega$ for each scenario, as in
previous reformulations, we obtain the master problem in iteration $k$ as
\begin{equation}\label{eq:2SSMILP-MP}
\begin{aligned}
	\min & \; c^\top x + \sum_{\omega \in \Omega} p_\omega z_\omega \\
  	\text{s.t.} &\ A^1 x \geq b^1 \;\\
        & \; z_\omega \geq \max_{1 \leq i \leq
	k} {\uphiIP^i}_{\omega}(b^2_\omega - A^2_\omega x)
        \quad \forall \omega \in \Omega \\ 
	&\ x \in X.
\end{aligned}
\end{equation}

Because each scenario dual function is the minimum of a collection of affine
functions, the overall master problem can be eventually reformulated as an
MILP by introducing additional binary variables (see~\cite{HasRal14a} for
details).

\subsubsection{Overall Algorithm \label{subsubsec:2SSMILPsAlgo}}

Putting this all together, in each iteration $k$, a master problem of the
form~\eqref{eq:2SSMILP-MP} is solved to obtain its optimal solution $(x^k,
\{z^k_\omega\}_{\omega \in \Omega})$ and a lower bound. Following that,
the subproblem is solved, which consists of evaluating $\phiIP(b^2_\omega -
A^2_\omega x^k)$ for each $\omega \in \Omega$ using a branch-and-bound
algorithm. The result is a strong
dual function~\eqref{eq:MILP-DF} for each scenario, as well as an overall
upper bound. If the upper and lower bounds are equal, then we are done.
Otherwise, the dual functions are fed back into the master problem and the
method is iterated until the upper and lower bounds converge.

\section{Mixed Integer Bilevel Linear Optimization Problems \label{sec:MIBLPs}}

We now move on to a detailed discussion of the application of this
generalization of Benders' principle to MIBLPs. As described earlier, MIBLPs
are two-stage MMILPs in which the variables at each stage are conceptually
controlled by different DMs with different objective functions.
MIBLPs model problems in game theory, specifically
the \emph{Stackelberg games} introduced by~\cite{stackelberg34}. Bilevel
optimization problems in the form presented here were formally introduced and
the term was coined in the 1970s by
\cite{bracken-mcgill73}, but computational aspects of such
optimization problems have been studied since at least the 1960s (see,
e.g.,~\cite{wollmer64}). Most of the initial research was limited to
continuous bilevel linear optimization problems containing only continuous
variables and linear constraints in both the stages.

Study of bilevel optimization problems containing integer variables and
algorithms for solving them is generally acknowledged to have been initiated
by~\cite{moorebard90}, who discussed the computational challenges of solving
such problems and suggested one of the earliest algorithms, a branch-and-bound
algorithm, which converges to an optimal solution if all first-stage variables
are integer or all second-stage variables are continuous.
Since then, many works have focused on special cases, such as those in which
the first-stage variables are all binary or all second-stage variables are
continuous. It is only in the past decade or so that study of exact algorithms
for the general case in which there are both continuous and general integer
variables in both stages has been undertaken.

Table~\ref{table:previousWork} provides a timeline of the main developments in
the evolution of such exact algorithms, indicating the types of
variables supported in both the first and second stages (C indicates
continuous, B indicates binary, and G indicates general integer). Most of
these works are either pure cutting plane or branch-and-cut algorithms in the
full space of first- and second-stage variables, and hence, are not technically
included under the umbrella of the framework of this paper.
Only three works, \cite{saharidis-ierapetritou08}, \cite{zengan14} and
\cite{yueetal19}, present themselves as decomposition algorithms. Of these
three, only the first work is a pure Benders-type algorithm, but it focuses on
the special case with all continuous second-stage variables,
in which case the reformulation can be done using
standard KKT conditions and the Benders' cuts are linear. The other two works
deviate from our approach in significant ways since their master problems
are in the full space of first- and second-stage variables.

Although no existing algorithm for the general MIBLP case can be considered as
a pure Benders-type algorithm, there nevertheless must necessarily
be \emph{some} connection between all algorithms for solving MIBLPs because of
the necessity to at least implicitly construct \emph{primal} approximations of
the MILP value function~\eqref{eq:MILP-VF}, a topic introduced in
Section~\ref{sec:reaction-fcn}. In the particular case when the
objective functions for the two stages disagree, such primal approximation of
the value function cannot be avoided. But while this need for primal
approximation may make it \emph{appear} as if some algorithmic alternatives
are also in fact Benders-type algorithms, it is the need for
explicit \emph{dual} approximations that sets such algorithms apart. The dual
approximation is necessary precisely \emph{because} of the projection
operation that is necessary when the second-stage variables are not present in
the master problem. Once second-stage variables are present in the master
problem, the dual approximation is no longer needed. With respect to the specific
algorithmic step of constructing primal functions, the construction methods of
existing works can be seen as special cases of our method, which is to
construct a parametric primal function~\eqref{eq:milp-ubf} (see
Section~\ref{sec:reaction-fcn}). For example, the primal functions in \cite{lozanosmith17}
and \cite{yueetal19} are (non-parametric) constant functions, and those
in \cite{capraraetal16} are a special case that exploits the specific structure
of the problem to obtain parametric functions that are linear in the first-stage
variables.

\begin{table}
\centering
\caption[Previous Work]{Evolution of algorithms for bilevel
  optimization \label{table:previousWork}}
\begin{tabular}{lll}
\hline\noalign{\smallskip}
Citation & Stage 1 Variable Types & Stage 2 Variable Types \\
\noalign{\smallskip}\hline\noalign{\smallskip}
\cite{wen90} & B & C \\
\cite{bard92} & B & B\\
\cite{faisca-etal07} & B, C & B, C \\
\cite{garcesetal09} & B & C \\
\cite{saharidis-ierapetritou08} & B, C & C \\
\cite{DeNRal09}, \cite{DeNegre2011} & G & G \\
\cite{koppe10} & G or C & G \\
\cite{baringo-conejo12} & B, C & C \\
\cite{xuwang14} & G & G, C \\
\cite{zengan14} & G, C & G, C \\
\cite{caramiamari15} & G & G \\
\cite{capraraetal16} & B & B \\
\cite{HemSmi16} & B, C & B, C \\
\cite{lozanosmith17} & G & G, C \\
\cite{wangxu17} & G & G \\
\cite{fischettietal17b}, \cite{fischettietal17a} & G, C & G, C \\
\cite{yueetal19} & G, C & G, C \\
\cite{TahRalDeN20} & G, C & G, C \\
\noalign{\smallskip}\hline
\end{tabular}
\end{table}

\subsection{Formulation \label{sec:MIBLP-formulation}}
To state the class of problems formally, we must introduce a type of
constraint that cannot be expressed in the canonical language of 
mathematical optimization. In addition to the usual linear
constraints, we have a constraint that requires the second-stage solution
to be optimal with respect to a problem that is parametric in the first-stage
solution. The formulation including this constraint, as it usually appears in
the literature on bilevel optimization, is
\begin{equation}\label{eq:miblp}
\begin{aligned}
	\min & \; c^\top x + {d^1}^\top y\\
        \text{s.t.} & \text{ } A^1x + G^1 y \geq b^1 \\
        & \; x \in X \\
	& y \in \arg\min \;\{{d^2}^\top \check{y} \\
        & \quad \text{s.t.} \; G^2\check{y} \geq b^2 -A^2 x \\
        & \quad \check{y} \in Y\},
\end{aligned}
\end{equation}
where
$A^1 \in \Q^{m_1 \times n_1}$, $G^1 \in \Q^{m_1 \times n_2}$, $b^1 \in
\Q^{m_1}$, $A^2 \in \Q^{m_2 \times n_1}$, $G^2 \in \Q^{m_2 \times n_2}$, $b^2
\in \Q^{m_2}$, $c \in \Q^{n_1}$, $d^1 \in \Q^{n_2}$, and $d^2 \in \Q^{n_2}$.
Note that the above-mentioned parametric problem
is nothing but the evaluation of $\phiIP(b^2 - A^2 x)$, where $\phiIP$ is the MILP
value function~\eqref{eq:MILP-VF}.

Underlying the above formulation are a number of assumptions. First, there is
an implicit assumption that whenever the evaluation of $\phiIP(b^2 - A^2 x)$
yields multiple optimal solutions, the one that is most advantageous for the
first-stage DM is chosen. This form of MIBLP is known as the \emph{optimistic}
case and is just one of several variants. The \emph{pessimistic} variant, for
example, is one in which the second-stage solution chosen is always the one
\emph{least} advantageous for the first-stage DM. It should also be pointed
out that we explicitly allow the second-stage variables in the constraints
$A^1 x + G^1 y \geq b^1$. This is rather non-intuitive but there are
applications for which this is a necessary element. We now define so-called
\emph{linking variables}.

\begin{definition}[\textbf{Linking Variables}]
Linking variables are the first-stage variables whose indices are in the set
\begin{equation*}
 L = \left\{i \in \{1, \hdots, n_1\} \mid A^2_i \neq 0\right\},
\end{equation*}
where $A^2_i$ denotes the $i^\textrm{th}$ column of $A^2$.
\end{definition}

\begin{assumption}\label{as:existence}
All linking variables are integer variables.
\end{assumption}

This assumption is required to ensure the existence of an optimal solution for
the given MIBLP. The optimal solution may not be attainable when there
are linking variables that are continuous and second-stage variables that are
integer~\citep{moorebard90, vicente96, koppe10}.

\begin{assumption}\label{as:equivalence}
All first-stage variables are linking variables.
\end{assumption}

Since we focus on optimistic bilevel problems, all non-linking variables can
simply be moved to the second stage
without loss of generality. This
assumption is made primarily for ease of
exposition, nevertheless results in a mathematically equivalent MIBLP, despite
the inconsistency with the intent of the original model. In combination
with Assumption~\ref{as:existence}, this assumption implies that all
first-stage variables are integer variables, i.e., $n_1 = r_1$.

\begin{assumption}\label{as:FRbounded}
The set
\begin{equation*}
\{(x, y) \in \Re^{n_1}_+ \times \Re^{n_2}_+ \mid y \in \Pcal_1(b^1 - A^1 x) \cap
\Pcal_2(b^2 - A^2 x)\} 
\end{equation*}
\end{assumption}
is bounded, where
\begin{equation*}
\Pcal_1(\beta^1) = \left\{y \in \Re_+^{n_2} \midd G^1 y \geq \beta^1  \right\}
\end{equation*}
and $\Pcal_2(\beta^2)$ (as defined in~\eqref{eq:MILP-VF-FR})
represent families of polyhedra consisting of all points satisfying
$G^1 y \geq \beta^1$ and $G^2 y \geq \beta^2$ for given right-hand sides
$\beta^1 \in \Re^{m_1}$ and $\beta^2 \in \Re^{m_2}$.
Assumption~\ref{as:FRbounded} is made to avoid uninteresting cases involving
unboundedness, but is easy to relax in practice.

\begin{assumption}\label{as:L2bounded}
For all $x \in \Re^{n_1}$, we have
\begin{equation*}
\phiIP(b^2 - A^2 x) > -\infty,
\end{equation*}
or, equivalently
\begin{equation*}
\left\{r \in \Re^{n_2}_+ \midd G^2 r \geq 0, {d^2}^\top r < 0\right\} = \emptyset.
\end{equation*}
\end{assumption}

Assumption~\ref{as:L2bounded} is also made to avoid uninteresting cases
involving unboundedness. Observe that in the case $\phiIP(b^2 - A^2 x) = -\infty$
for a given value of $x$, then $\phiIP(b^2 - A^2 x) = -\infty$ for all values
of $x$.
Note that Assumptions~\ref{as:equivalence}-\ref{as:L2bounded} can be relaxed
in practice.

\subsection{Projection \label{sec:reaction-fcn}}
We now apply the familiar operations of Benders' framework to the formulation
\eqref{eq:miblp}. Upon projecting into the space of the first-stage  variables,
we obtain the reformulation
\begin{equation}\label{eq:miblp-ref}
\min_{x \in X} \left\{c^\top x + \rho(b^1 - A^1x,\;b^2 - A^2 x) \right\},
\end{equation}
where $\rho$ is the second-stage \emph{reaction function}, defined as
\begin{equation}\label{eq:ReF}
\rho(\beta^1, \beta^2) = \min \left\{{d^1}^\top y \midd y \in \Pcal_1(\beta^1),\;y
\in \arg\min \left\{{d^2}^\top \check{y} \midd \check{y} \in \Pcal_2(\beta^2) \cap Y \right\}
\right\}.
\end{equation}

Observe that the reaction function
has embedded within it exactly the kind of optimality constraint we tried
to eliminate by using projection to reformulate the given MIBLP.
Although the reaction function appears at first to be the value function of a
general bilevel optimization problem, it is actually the value function of a
lexicographic MILP. To see this, note that the evaluation of
$\rho(\beta^1, \beta^2)$ for known values of $(\beta^1, \beta^2)$ can be done
in two steps. First, we evaluate
MILP
\begin{equation*}
\phiIP(\beta^2)
= \min \left\{{d^2}^\top \check{y} \midd \check{y} \in \Pcal_2(\beta^2) \cap
Y \right\}. 
\end{equation*}
Then, we have
\begin{equation}\label{eq:ReF-VF}
\rho(\beta^1, \beta^2) = \min \left\{{d^1}^\top y \midd
{d^2}^\top y \leq \phiIP(\beta^2),
y \in \Pcal_1(\beta^1) \cap \Pcal_2(\beta^2) \cap Y \right\},
\end{equation}
which is an MILP. In other words, the evaluation of $\rho(\beta^1, \beta^2)$
is an optimization problem over the set of optimal solutions to a given
(non-parametric) MILP. Once $\phiIP(\beta^2)$ is known, the set of optimal
solutions over which we are trying to optimize is the feasible set of an MILP.
The reason this is not a bilevel optimization problem is simply because the
right-hand side vector $\beta^2$ is not parametric, i.e., it is a known vector.
In the next part of this section, we examine the properties and structure of
the reaction function before discussing how to construct associated dual
functions.

\paragraph{Reaction Function.} As with all value functions, $\rho(\beta^1,
\beta^2) = \infty$ for a given $(\beta^1, \beta^2) \in \Re^{m_1} \times
\Re^{m_2}$ if either $\{y \in \Re_+^{n_2} \mid y \in \Pcal_1(\beta^1) \cap
\Pcal_2(\beta^2) \cap Y\} = \emptyset$ or $\phiIP(\beta^2) = -\infty$ (which
cannot happen under Assumption~\ref{as:L2bounded}), and $\rho(\beta^1,
\beta^2) = -\infty$ for all $(\beta^1, \beta^2) \in \Re^{m_1} \times
\Re^{m_2}$ if the lexicographic MILP is itself unbounded, i.e., $\{r \in
\Re^{n_2}_+ \mid G^1 r \geq 0, G^2 r \geq 0, {d^1}^\top r < 0\} \neq \emptyset$
(which cannot happen under Assumption~\ref{as:FRbounded}).

We now illustrate the structure of the reaction function with a simple
example. Although its structure is combinatorially more complex than that of
the MILP value function, it nevertheless also has a piecewise polyhedral
structure. We do not provide formal results concerning the structure and
properties of the reaction function here, but these can be derived by
application of techniques similar to those used to derive the structure of the
MILP value function.
\begin{example}
Consider the following reaction function arising from an MIBLP with $G^1 = 0$.
\begin{equation}\label{eq:toyReF}
\begin{aligned}
	\rho(\beta) &=& \min& \;  - y_1 + y_2 - 5 y_3 + y_4\\
         && \text{s.t.}\; & (y_1, y_2, y_3, y_4) \in \arg\min \;\{2 \check y_1 +
         4 \check y_2 + 3 \check y_3 + 4 \check y_4\\
         &&&\qquad\qquad\qquad\qquad \text{s.t.}\; 2 \check y_1 + 5 \check y_2 +
         2 \check y_3 + 2 \check y_4 \geq \beta \\
         &&& \qquad\qquad\qquad\qquad\qquad \check y_1, \check y_2, \check y_3
         \in \Z_+,\ \check y_4 \in \Re_+\}
\end{aligned}
\end{equation}

This function can be reformulated as
\begin{equation*}
\begin{aligned}
	\rho(\beta) &=& \min& \; - y_1 + y_2 - 5 y_3 + y_4\\
         && \text{s.t.}\; & 2 y_1 + 5 y_2 + 2 y_3 + 2 y_4 \geq \beta \\
         &&& 2 y_1 + 4 y_2 + 3 y_3 + 4 y_4 \leq \phiIP(\beta)\\
         &&& y_1, y_2, y_3 \in \Z_+,\ y_4  \in \Re_+.
\end{aligned}
\end{equation*}
using the MILP value function $\phiIP$, which is the same as for~\eqref{eq:ip}.
The function $\rho$ is plotted in Figure~\ref{fig:ToyReF}.
\begin{figure}
\includegraphics[width=\textwidth]{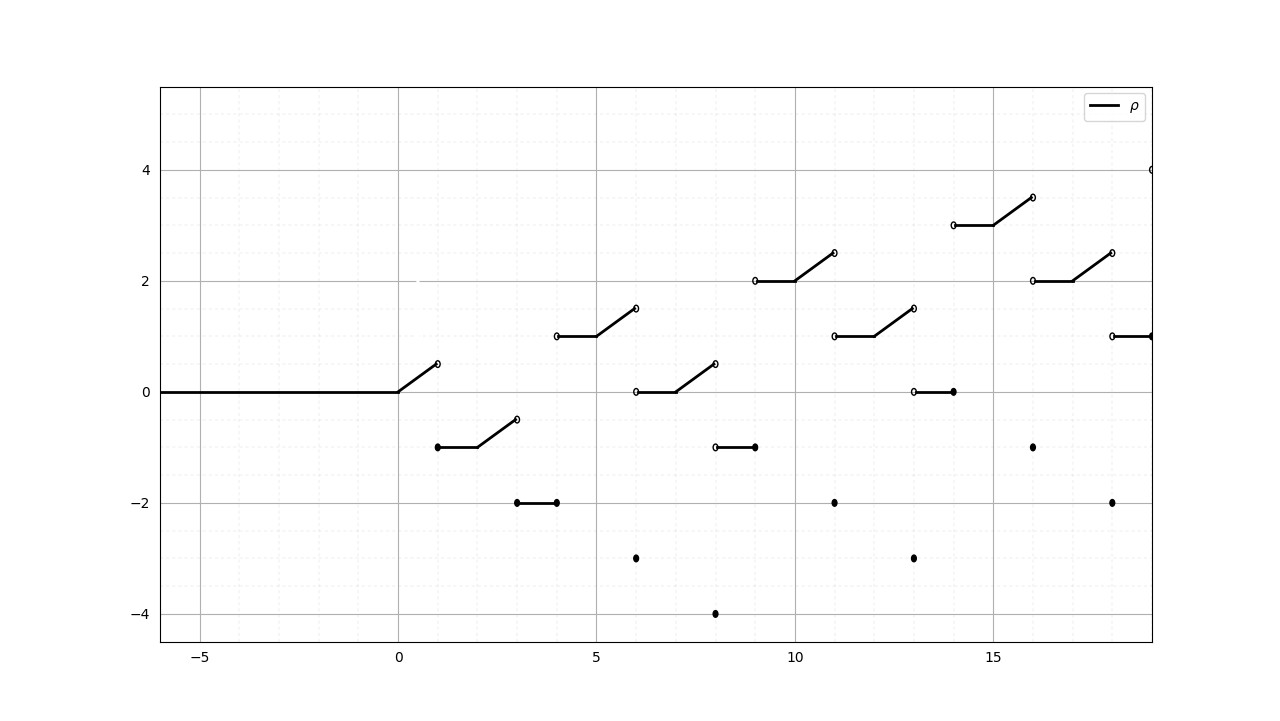}
\caption{Reaction function $\rho$ \eqref{eq:toyReF} \label{fig:ToyReF}}
\end{figure}
\qed
\end{example}

As usual, we do not have an exact description of $\rho$ in general, so
we cannot solve~\eqref{eq:miblp-ref} directly, but instead replace $\rho$
in~\eqref{eq:miblp-ref} with a dual function $\urho$. Following the
earlier procedure, this dual function is taken to be the maximum of the strong
dual functions $\urho^k$ obtained by solving a subproblem in each
iteration $k$. The resulting master problem in iteration $k$ is
\begin{equation}\label{eq:MIBLP-MP}
\begin{aligned}
	& & \min & \; c^\top x + z \\
        & & \text{s.t.} & \; z \geq \urho^i(b^1 - A^1 x,\;b^2 - A^2 x), \quad 1 \leq i \leq k\\
	& & & x \in X.
\end{aligned}
\end{equation}

Similarly, the subproblem in iteration $k$ is to evaluate $\rho(b^1 - A^1
x^k,\;b^2 - A^2 x^k)$ for the solution $x^k$ to~\eqref{eq:MIBLP-MP}, in order
to construct a dual function $\urho^k$ that is strong at $(b^1 - A^1
x^k,\;b^2 - A^2 x^k)$. We next detail the construction of this strong dual
function.

\paragraph{Dual Functions. \label{sec:reaction-fcn-dual}}
As we have already noted, the subproblem in iteration $k$ is to evaluate the
reaction function~\eqref{eq:ReF} for $(\beta^1, \beta^2) = (b^1 - A^1 x^k, b^2
- A^2 x^k)$. This problem is an MILP and we have the
following theorem based on Theorem~\ref{thm:milp-lbf}.
\begin{theorem} \label{thm:ReFLowerApprox}
Let $(\hat{\beta^1}, \hat \beta^2) \in \Re^{m_1} \times \Re^{m_2}$ be such
that $\rho(\hat{\beta^1}, \hat \beta^2) < \infty$ and suppose $T$ is the set
of indices of leaf nodes of a
branch-and-bound tree resulting from evaluation of $\rho
(\hat{\beta^1}, \hat \beta^2)$. Then there exists a dual function
$\utilde{\rho} : \Re^{m_1} \times \Re^{m_2} \rightarrow \Re \cup \{\pm
\infty\}$ of the form
\begin{equation} \label{eq:ref-lbf}
\utilde{\rho} (\beta^1, \beta^2) = \min_{t \in T} \left({\beta^1}^\top
          \eta^{1,t} + {\beta^2}^\top \eta^{2,t} + \phiIP(\beta^2) \eta^{\phi,t} +
          \alpha^t \right) \quad \forall 
(\beta^1, \beta^2) \in \Re^{m_1} \times \Re^{m_2},
\end{equation}
where $(\eta^{1,t}, \eta^{2,t}, \eta^{\phi,t}) \in \Re^{m_1} \times \Re^{m_2}
\times \Re$ is a dual feasible solution of the LP relaxation associated with
node $t$, and $\alpha^t \in \Re$ is the product of reduced costs and variable
bounds of this LP relaxation. Further, this dual function is strong at $(\hat
\beta^1, \hat \beta^2)$ if $\utilde{\rho} (\hat{\beta^1}, \hat \beta^2)=
\rho(\hat{\beta^1}, \hat \beta^2)$.
\end{theorem}

The interpretation of this result is similar to the interpretation of
Theorem~\ref{thm:milp-lbf}. Let us look at an example that depicts these
functions.

\begin{example} \label{eg:reaction-fcn-dual}
Figure~\ref{fig:ToyReF2LBFs} shows five dual functions obtained
upon applying the result in Theorem~\ref{thm:ReFLowerApprox} to the reaction
function~\eqref{eq:toyReF} ($\utilde{\rho}^1$ for $\beta = 0$,
$\utilde{\rho}^2$ for $\beta = 1$, $\utilde{\rho}^3$ for $\beta = 2$,
$\utilde{\rho}^4$ for $\beta = 5$, and $\utilde{\rho}^5$ for $\beta = 8$).
\begin{figure}
\includegraphics[width=\textwidth]{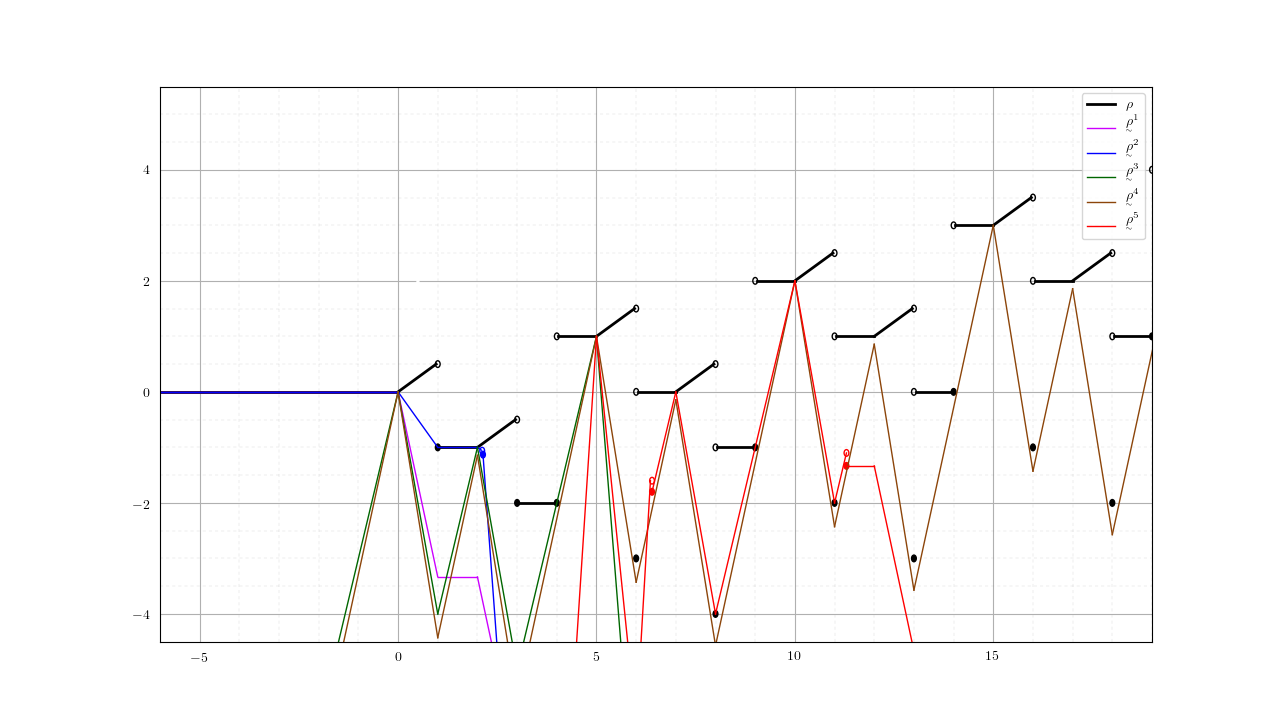}
\caption{Dual functions for \eqref{eq:toyReF} \label{fig:ToyReF2LBFs}}
\end{figure}
As expected, these dual functions are piecewise polyhedral. For example,
solving~\eqref{eq:toyReF} with $\beta = 8$ as an equivalent MILP (after
obtaining $\phiIP(8)$ at first) yields the dual information (dual
solution and reduced costs) from leaf nodes of the branch-and-bound tree shown
in Table~\ref{tab:reaction-fcn-dual}.
\begin{table}
\centering
\caption{Data for construction of the dual function in
  Example~\ref{eg:reaction-fcn-dual} \label{tab:reaction-fcn-dual}}
\begin{tabular}{cccccc}
\hline\noalign{\smallskip}
$t$ & Branching constraint & $(\eta^{2,t}, \eta^{\phi,t})$ & $\underline\eta^t$
  & $\overline\eta^t$ & ${\eta^{2,t}}^\top\beta +  \eta^{\phi,t}\phiIP(\beta) +
  \alpha^t$\\
\noalign{\smallskip}\hline\noalign{\smallskip}
1 & $y_2 \geq 2$ & $(0, -\frac{5}{3})$ & $(\frac{7}{3}, \frac{23}{3}, 0,\
  \frac{7}{3})$ & $(0, 0, 0, 0)$ & $-\frac{5}{3}\phiIP(\beta) + \frac{46}{3}$ \smallskip\\
2 & $y_2 \leq 1$, $y_1 \leq 0$ & $(11, -9)$ & (0, 0, 0, 15) & (-5, -18, 0, 0) &\
  $11 \beta - 9 \phiIP(\beta) - 18$ \smallskip\\
3 & $y_2 \leq 1$, $y_1 \geq 1$, $y_3 \leq 0$ & $(3, -3.5)$ & (0, 0, 0, 9) & (0,\
  0, -0.5, 0) & $3 \beta - 3.5 \phiIP(\beta)$ \smallskip\\
4 & $y_2 \leq 1$, $y_1 \geq 1$, $y_3 \geq 1$ & $(13, -16)$ & (5, 0, 17, 39) &\
  (0, 0, 0, 0) & $13 \beta - 16 \phiIP(\beta) + 22$ \\
\noalign{\smallskip}\hline
\end{tabular}
\end{table}

This results in the dual function
\begin{equation*}
\utilde{\rho}^5(\beta) = \min\left\{-\frac{5}{3}\phiIP(\beta) + \frac{46}{3}, 11\
  \beta - 9 \phiIP(\beta) - 18, 3 \beta - 3.5 \phiIP(\beta), 13 \beta - 16\
  \phiIP(\beta) + 22\right\},
\end{equation*}
containing the MILP value function $\phiIP$. The remaining dual functions are
obtained the same way.
\qed
\end{example}

It is clear from Theorem~\ref{thm:ReFLowerApprox} that the construction of
$\utilde{\rho}$ in~\eqref{eq:ref-lbf} implicitly requires construction of the
value function $\phiIP$. However, the construction of $\phiIP$ is itself a
difficult task and generally impractical. Further, the complex structure of
$\phiIP$ makes the structure of $\utilde{\rho}$ highly complex. To work around
this difficulty, we replace $\phiIP$ in~\eqref{eq:ReF-VF} with
a \emph{primal function}, which bounds the value
function from above and is strong at the given right-hand side. This
replacement results in an alternative dual function (which we denote by
$\urho$) that is still strong at the given right-hand side. We use
$\urho$ in place of $\utilde{\rho}$ in our work. To this end, we
now embark on a small diversion into MILP primal functions.

\paragraph{Primal Functions. \label{subsubsec:MILPUpperApprox}}

In contrast with dual functions, strong \emph{primal functions} bound the
value function from above.
\begin{definition}[\textbf{Primal Function}]\label{def:primalfn}
A primal function $P : \Re^{m_2} \rightarrow \Re \cup \{\pm \infty\}$ is one
that  
satisfies $P(\beta) \geq \phiIP(\beta)$ for all $\beta \in \Re^{m_2}$.  It is
\emph{strong} at $\hat \beta \in \Re^{m_2}$ if $P(\hat \beta) =
\phiIP(\hat \beta)$.
\end{definition}

An obvious way to construct such a function is to consider the value function
of a restriction of the given MILP (see~\cite{Guzelsoy2009} and
\cite{GuzRal07} for methods of construction). The following theorem
presents the main result for constructing strong primal functions 
from restrictions of the given MILP.

\begin{theorem} [\cite{Guzelsoy2009}, Theorem 3.39] \label{thm:milp-ubf}
Consider the MILP value function~\eqref{eq:MILP-VF}. Let $K \subseteq N
\coloneqq \{1, \hdots, n_2\}, \;s \in \Re^{|K|}_+$ be
given, and define the function $\bphiNK : \Re^{m_2} \ra \Re \cup \{\pm
\infty\}$ such that 
\begin{equation*}
\bphiNK(\beta) = \sum_{i \in K} d^2_i s_i + \phi_{N\sm K}\left(\beta - \sum_{i \in
  K}G^2_i s_i\right) \quad \forall \beta \in \Re^{m_2},
\end{equation*}
where $G^2_i$ is the $i^{th}$ column of $G^2$ and
\begin{equation*}
\begin{aligned}
\phi_{N\sm K}(\beta) & = & \min & \; \sum_{i \in N\sm K} d^2_i y_i \\
        & & \text{s.t.} & \; \sum_{i \in N\sm K}G^2_i y_i \geq \beta \\
	& & &\ y_i \in \Z_+ \text{ } \forall i \in I, \text{  } y_i \in
\Re_+ \text{ } \forall i \in C, 
\end{aligned}
\end{equation*}
where $I \subseteq N$ and $C \subseteq N$ represent sets of indices for integer
and continuous variables respectively. Then, $\bphiNK$ is a valid primal
function of $\phiIP$, i.e., $\bphiNK(\beta) \geq \phiIP(\beta) \text{ } \forall
\beta \in \Re^{m_2}$, if $s_i \in \Z_+ \text{ } \forall i \in I \cap K$ and
$s_i \in \Re_+ \text{ } \forall i \in C \cap K$. Further, the function $\bphiNK$
will be strong at a given right-hand side $\hat \beta \in \Re^{m_2}$ if and only
if $s_i = y_i^* \; \forall i \in K$ where $y^*$ is an optimal solution of $\phiIP(\hat \beta)$.
\end{theorem}

By convention, for a known $\hat\beta \in \Re^{m_2}$, we consider
$\phi_{N\sm K}(\hat\beta) = \infty$ if the corresponding problem is infeasible
and $\phi_{N\sm K}(\hat\beta) = - \infty$ if the problem is unbounded.
The above result indicates that a primal function obtained from a restriction
in which the values of certain variables have been fixed is strong at
$\hat{\beta} \in \Re^{m_2}$ if the fixed values of these variables
correspond to those of an optimal solution at $\hat{\beta}$.
A convenient approach is to fix all the integer variables to their optimal
values. If there are no continuous variables in the problem, then the resulting
primal function is
\begin{equation*}
\bar{\phi}_{\emptyset}(\beta) = \left\{\begin{array}{ll}
\sum\limits_{i \in I} d^2_i y_i^* & \text{if}\;\; \beta = \hat\beta\\
\infty & \text{otherwise}
\end{array}\right.,
\end{equation*}
which is a single point, but still a valid strong primal function at
$\hat\beta$.  If continuous variables exist, then the restriction is a
continuous restriction mentioned earlier. The resulting value function
$\phi_{C}$ is nothing but the value function of an LP discussed briefly in
Section~\ref{subsec:LPIllustration}.
Let us now look at an example of using continuous restrictions to generate
primal functions.
\begin{example}\label{eg:toyMILPUBFs}
Consider the MILP~\eqref{eq:ip}. Figure~\ref{fig:IPVF2UBFs} demonstrates four
primal functions obtained upon applying the result in
Theorem~\ref{thm:milp-ubf} to this MILP. Specifically, $\bphi{1}$ for
$\beta = 0$, $\bphi{2}$ for $\beta = 2$, $\bphi{3}$ for $\beta = 4$,
and $\bphi{4}$ for $\beta = 5$. Here, we
consider continuous restrictions of the given MILP. For example, for $\beta =
0$, the optimal solution of~\eqref{eq:ip} is $(0, 0, 0, 0)$ resulting in the
following continuous restriction.
\begin{equation*}
\phi_{C}(\beta) = \min \left\{4 y_4 \midd 2 y_4 \geq \beta, \;y_4 \in \Re_+
  \right\} \quad \forall \beta \in \Re
\end{equation*}

It is easy to observe that this LP value function is
\begin{equation*}
\phi_{C}(\beta) = \left\{\begin{array}{ll}
0 & \text{if}\;\; \beta \leq 0\\
2 \beta & \text{otherwise}
\end{array}\right.,
\end{equation*}
which itself is the required primal function $\bphi{1}$ because
the integer component of the MILP objective function is zero. Other primal
functions can be constructed in a similar way by solving the MILP with a new
right-hand side, calculating the integer component of the MILP objective
function at its optimal solution, and simply translating $\phi_{C}$ based on
this integer component value.
\begin{figure}
\includegraphics[width=\textwidth]{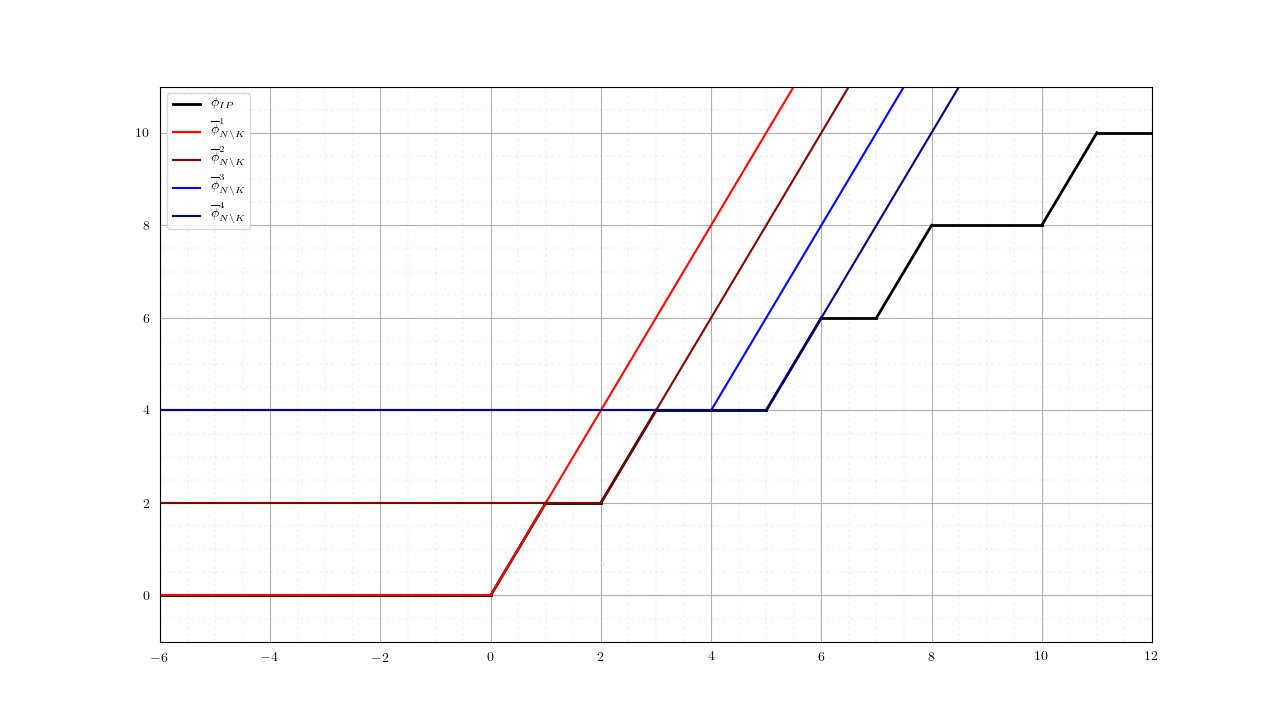}
\caption{Primal functions for~\eqref{eq:ip} \label{fig:IPVF2UBFs}}
\end{figure}
The primal function can be strengthened in the same way as dual functions are
strengthened, by considering the minimum of multiple such (strong) functions.
The epigraph of such a function is the minimum of convex radial cones and
equals the epigraph of the value function when enough such cones are
considered, as mentioned in Section~\ref{subsubsec:2SSMILPs-proj}.
\qed
\end{example}

Although LPs are themselves easy to solve, obtaining a full description of the
LP value function $\phi_C$ is still difficult, so we consider a partial
function that we can easily obtain. Specifically, we use a single
hyperplane of~\eqref{eq:LP-VF-Struct} to construct $\bphiNK$. In
Example~\ref{eg:toyMILPUBFs}, this is equivalent to considering only one of
the two hyperplanes forming the cone corresponding to $\phi_C$. It is obvious
from Figure~\ref{fig:IPVF2UBFs} that any single hyperplane cannot form a valid
primal function in the entire domain of the right-hand side vector. Therefore,
we also need to restrict the domain of the right-hand side vector to only the
region in which the single hyperplane is a valid upper bound.

Let $y^* = (y^*_I, y^*_C)$ be an optimal solution of an instance
of~\eqref{eq:ReF} for a known right-hand side $(\hat\beta^1, \hat\beta^2)$,
where $I$ and $C$ correspond to sets of indices of second-stage integer and
continuous variables respectively. This inherently implies that
$\phiIP(\hat\beta^2) = {d^2}^\top y^*$. The value function $\phi_C$ of the
continuous restriction is then
\begin{equation*}
\begin{aligned}
	& & \phi_C(\beta) = \min & \; {d^2_C}^\top y_C \\
        & & \text{s.t.} & \; G^2_C y_C \geq \beta \\
	& & &\ y_C \geq 0.
\end{aligned}
\end{equation*}

Let $\eta^*$ be an optimal solution of the dual of this LP, with right-hand
side $\beta^2 - G^2_I y^*_I$. Then ${\eta^*}^\top \beta$ is a dual function
strong at $\beta^2 - G^2_I y^*_I$. From the theory of LP duality, we know that
this function provides a valid upper bound as long as $\eta^*$ remains
optimal, which is the case for all $\beta$ such that $(G^2_B)^{-1} \beta \geq
0$, where $B$ is the index set corresponding to the optimal basis and
${G^2_B}$ is the optimal basis matrix.

Thus, we obtain our final primal function with a restricted domain
\begin{equation}\label{eq:milp-ubf}
\bphiIP(\beta^2) = \left\{\begin{array}{ll}
\left(\beta^2 - G^2_I y^*_I\right)^\top{\eta^*} + {d^2_I}^\top y^*_I &
  \text{if}\; {\left(G_B^2\right)}^{-1}
\left(\beta^2 - G^2_I y^*_I\right) \geq 0\\
\infty & \text{otherwise}
\end{array}\right.
\end{equation}
with which we replace $\phiIP$ in~\eqref{eq:ref-lbf}.  This
ensures that the dual function $\urho$ that we construct in each
iteration of the algorithm is valid for all values of $(\beta^1, \beta^2)$.

\subsection{Master Problem}
Combining the results obtained in the previous section, the final form of
$\urho$ is
\begin{equation}\label{eq:ref-lbf-final}
\urho(\beta^1, \beta^2) = \min_{t \in
T} \left({\beta^1}^\top \eta^{1,t} 
+ {\beta^2}^\top \eta^{2,t} + \bphiIP(\beta^2) \eta^{\phi,t} + \alpha^t
\right),
\end{equation}
where $\bphiIP$ is the primal function~\eqref{eq:milp-ubf}. This results in the
optimality constraint
\begin{equation}\label{eq:MIBLP-OC}
z \geq \min_{t \in T} \left({\beta^1}^\top \eta^{1,t} +
{\beta^2}^\top \eta^{2,t} +
\bphiIP(\beta^2) \eta^{\phi,t} + \alpha^t
\right)
\end{equation}
that we add to the master problem in each iteration of the algorithm, with
$(\beta^1, \beta^2) = (b^1 - A^1 x, b^2 - A^2 x)$. Finally, the
updated master problem after iteration $k$ of the algorithm is
\begin{equation}\label{eq:MIBLP-MP-Updated}
\begin{aligned}
	& & \min & \; c^\top x + z \\
        & & \text{s.t.} & \; z \geq \min_{t \in T_i} \left\{ \left(b^1 - A^1
        x\right)^\top \eta^{1,t}_i 
  + \left(b^2 - A^2 x\right)^\top \eta^{2,t}_i  + \bphiIP^i(b^2 - A^2 x)
  \eta^{\phi,t}_i
        + {\alpha}^t_i\right\} \quad 1 \leq i \leq k\\ 
	& & &\ x \in X,
\end{aligned}
\end{equation}
where the vectors, matrices, sets, and functions with the subscript and
superscript $i$ correspond to the information obtained in iteration $i \leq k$
of the algorithm.

\subsection{Overall Algorithm \label{subsubsec:ourAlgo}}

We now have all the components required for solving~\eqref{eq:miblp} with the
generalized Benders' decomposition algorithm in Figure~\ref{fig:gpBenders}. In
each iteration $k$ of the algorithm, a master problem of the
form~\eqref{eq:MIBLP-MP-Updated} is solved to obtain its optimal solution $(x^k,
z^k)$ and a global lower bound. Then, the subproblem is solved as an equivalent
MILP, by evaluating~\eqref{eq:ReF} at $(b^1 - A^1 x^k, b^2
- A^2 x^k)$, to obtain a branch-and-bound tree and a global upper bound. Finally,
an optimality constraint of the form~\eqref{eq:MIBLP-OC} is constructed and
added to the master problem to strengthen $z$. This constraint introduces some
nonlinear components in the master problem but they can be linearized
(as mentioned below) to obtain an MILP formulation for the master problem. These
steps are repeated until the termination criterion is achieved.

We now illustrate the above discussion with an example.
\begin{example} \label{eg:miblp-toy}
Consider the MIBLP
\begin{equation}\label{eq:miblp-toy}
\begin{aligned}
	\min \;& x_1 - 3 x_2 - y_1 + y_2 - 5 y_3 + y_4\\
        \text{s.t.} \;& -x_1 + 2 x_2 \leq 1 \\
         & x_1 \leq 3,\; x_2 \leq 2,\; x_1, x_2 \in \Z_+ \\
   &(y_1, y_2, y_3, y_4) \in \arg\min \;\{2 \check y_1 + 4 \check y_2 + 3 \check
   y_3 + 4 \check y_4\\
         &\qquad\qquad\qquad\qquad \text{s.t.} \; 2 \check y_1 + 5 \check y_2 +
         2 \check y_3 + 2 \check y_4 \geq x_1 + x_2 \\
         &\qquad\qquad\qquad\qquad\qquad \check y_1, \check y_2, \check y_3, \in
         \Z_+,\ \check y_4 \in \Re_+\},
\end{aligned}
\end{equation}
which is based on~\eqref{eq:ip} and~\eqref{eq:toyReF}. Based on earlier
discussion, we solve four optimization problems in iteration $k$ of the
algorithm: a master problem, an MILP~\eqref{eq:ip} (with $\beta^k = x^k_1 +
x^k_2$), a subproblem~\eqref{eq:toyReF} (with $\beta^k = x^k_1 + x^k_2$), and a
continuous restriction of~\eqref{eq:ip}.

\textbf{Iteration 1.} Our initial dual function is simply
$\urho^0(\beta) = -\infty$ for all $\beta \in \Re^{m_2}$ and
solving the initial master problem yields the optimal solution $(x_1^1, x_2^1)
= (3, 2)$ and $z^1 = -\infty$, so that $\textrm{LB}^1 =-\infty$.
Then, we solve~\eqref{eq:ip} with right-hand side $x_1^1 + x_2^1 = 5$ to
obtain $\phiIP(x_1^1 + x_2^1) = 4$. Next, we solve the subproblem to obtain its
optimal solution $(y_1^1, y_2^1, y_3^1, y_4^1) = (0, 1, 0, 0)$, so we have
$\rho(x_1^1 + x_2^1) = 1$ and $\textrm{UB}^1 = x_1^1 - 3x_2^1 + \rho(x_1^1 +
x_2^1) = -2$. We obtain the dual information (dual solution, positive and
negative reduced costs) shown in Table~\ref{tab:miblp-toy} from the
branch-and-bound tree, which has only one node.
\begin{table}
\centering
\caption{Data for construction of the dual function in
  Example~\ref{eg:miblp-toy}. \label{tab:miblp-toy}}
\begin{tabular}{ccccc}
\hline\noalign{\smallskip}
$t$ & $(\eta_1^{2,t}, \eta_1^{\phi,t})$ & $(\underline\eta_1^{1,t},\
  \underline\eta_1^{2,t}, \underline\eta_1^{3,t}, \underline\eta_1^{4,t})$ &\
  $(\bar\eta_1^{1,t}, \bar\eta_1^{2,t}, \bar\eta_1^{3,t}, \bar\eta_1^{4,t})$ &\
  ${\eta_1^{2,t}}^\top \beta+ \eta_1^{\phi,t} \bphiIP^1(\beta) + \alpha_1^t$\\
\noalign{\smallskip}\hline\noalign{\smallskip}
1  & (3.29, -3.86) & $(0.14, 0, 0, 9.86)$ & $(0, 0, 0, 0)$ & $3.29 \beta - 3.86\
  \bphiIP^1(\beta)$ \\
\noalign{\smallskip}\hline
\end{tabular}
\end{table}

Since $\textrm{UB}^1 \not= \textrm{LB} ^1$, we further solve the continuous
restriction to obtain its optimal dual solution $\eta^*_1 = 0$ and optimal basis
inverse matrix $(G^2_{B,1})^{-1} = [-1]$. Finally, we construct and add the dual
function
\begin{equation*}
\urho^1(\beta) = \min\{3.29 \beta - 3.86 \bphiIP^1(\beta)\} = 3.29\
  \beta - 3.86 \bphiIP^1(\beta),
\end{equation*}
where
\begin{equation*}
\bphiIP^1(\beta) = \left\{\begin{array}{ll}
4 & \text{if}\;\; \beta \leq 5\\
\infty & \text{otherwise.}
\end{array}\right.,
\end{equation*}
to the master problem and proceed to the next iteration.

For conciseness, we now mention only $\urho^k(\beta)$ and
$\bphiIP^k(\beta)$ obtained in each iteration $k$.

\textbf{Iteration 2.}
\begin{equation*}
\urho^2(\beta) = - \frac{5}{3} \bphiIP^2(\beta)
\end{equation*}
\begin{equation*}
\bphiIP^2(\beta) = \left\{\begin{array}{ll}
0 & \text{if}\;\; \beta \leq 0\\
\infty & \text{otherwise}
\end{array}\right.
\end{equation*}

\textbf{Iteration 3.}
\begin{equation*}
\urho^3(\beta) = \min \left\{0.5 \bphiIP^3(\beta), -4.5\
 \bphiIP^3(\beta) + 8.5\right\}
\end{equation*}
\begin{equation*}
\bphiIP^3(\beta) = \left\{\begin{array}{ll}
2 & \text{if}\;\; \beta \leq 2\\
\infty & \text{otherwise}
\end{array}\right.
\end{equation*}

\textbf{Iteration 4.}
\begin{equation*}
\urho^4(\beta) = \min \left\{-\dfrac{5}{3} \bphiIP^4(\beta) +\
  \dfrac{23}{3}, 17 \beta - 13 \bphiIP^4(\beta),\ -0.5 \bphiIP^4(\beta),\
  21 \beta - 26 \bphiIP^4(\beta) + 40\right\}
\end{equation*}
\begin{equation*}
\bphiIP^4(\beta) = \left\{\begin{array}{ll}
4 & \text{if}\;\; \beta \leq 4\\
\infty & \text{otherwise}
\end{array}\right.
\end{equation*}

\textbf{Iteration 5.} Solving the updated master problem yields $(x^5_1, x^5_2,
z^5) = (1, 1, -1)$ and $\textrm{LB}^5 = -3$. Solving the subproblem yields $(y_1^5,
y_2^5, y_3^5, y_4^5) = (1, 0, 0, 0)$, $\rho(x_1^5 + x_2^5) = -1$ and $\textrm{UB}^5 =
-3$. Since $\textrm{UB}^5 = \textrm{LB}^5$, the termination criterion is
achieved. Hence, we stop the algorithm with an optimal solution $(x_1^*,
x_2^*, y_1^*, y_2^*, y_3^*, y_4^*) = (1, 1, 1, 0, 0, 0)$.

Figure~\ref{fig:IPVF4UBFsIter4} shows the value function $\phiIP$ and its primal
functions obtained in every iteration. Similarly,
Figure~\ref{fig:SPRF4LBFsIter4} shows the reaction function $\rho$ and its dual
functions. The function values are infinite wherever there is no plot. These
figures illustrate the fact that overall approximations of $\phiIP$ and $\rho$ are
strengthened after every iteration. After the final iteration, the
approximation function values are the same as the exact function values at the
right-hand side $x_1^* + x_2^* = 2$.
\begin{figure}
\includegraphics[width=\textwidth]{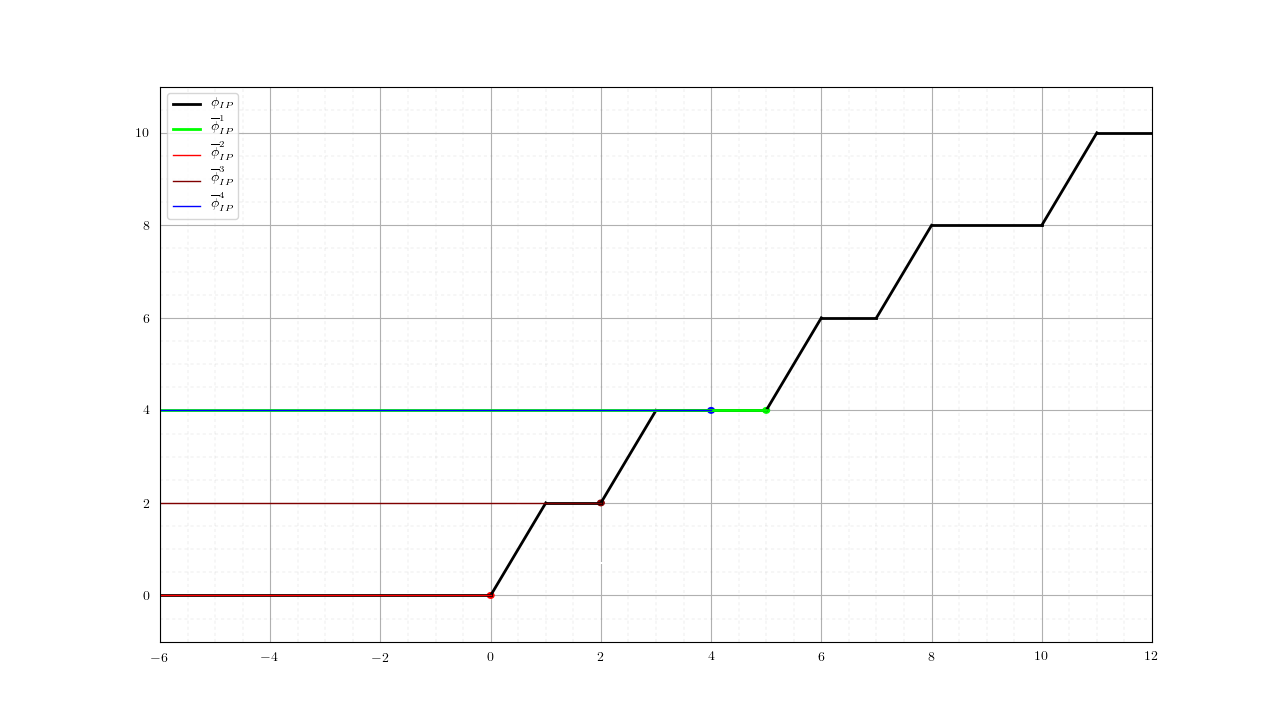}
\caption{Primal functions constructed in the algorithm for
  solving~\eqref{eq:miblp-toy} \label{fig:IPVF4UBFsIter4}}
\end{figure}
\begin{figure}
\includegraphics[width=\textwidth]{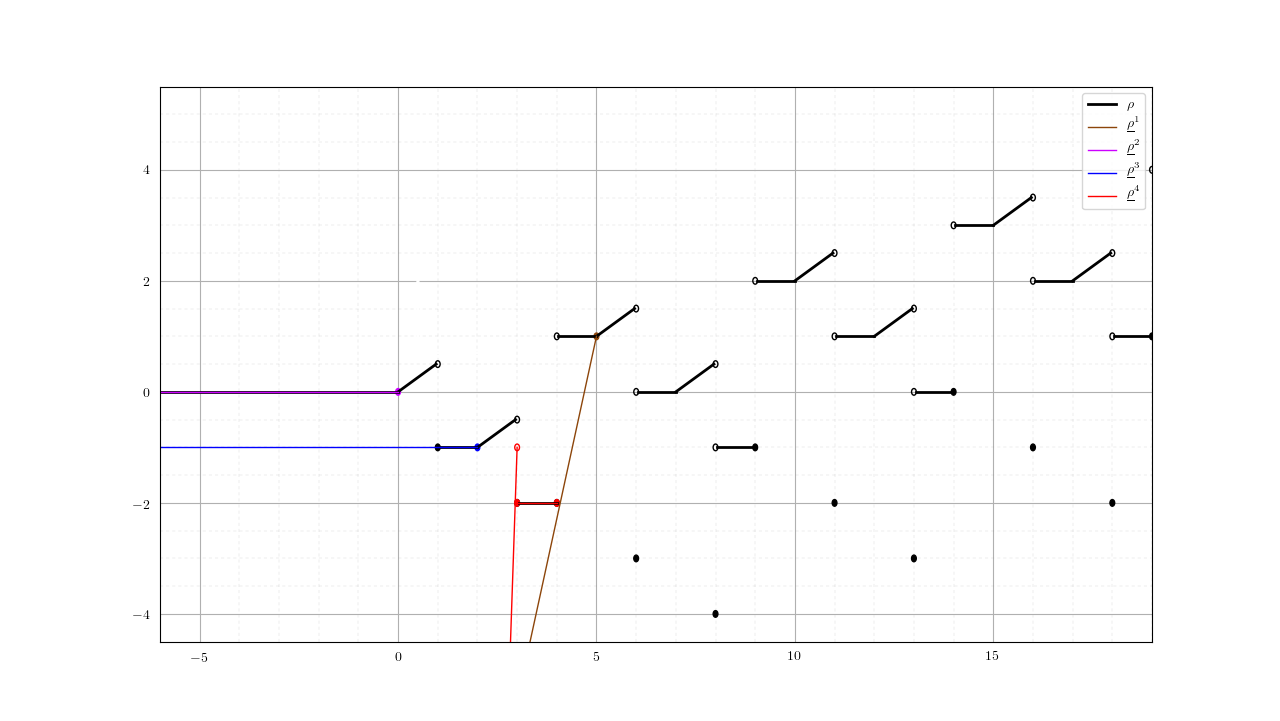}
\caption{Dual functions constructed in the algorithm for
  solving~\eqref{eq:miblp-toy} \label{fig:SPRF4LBFsIter4}}
\end{figure}
\qed
\end{example}

We now briefly discuss the linearization of the master problem
\eqref{eq:MIBLP-MP-Updated}. For notational simplicity, we drop the subscripts
and superscripts denoting the algorithmic iteration. There are two types of
nonlinearities in this problem: (1) the if-else condition
in~\eqref{eq:milp-ubf} and (2) the minimization operator
in~\eqref{eq:ref-lbf-final}. We eliminate these nonlinearities by introducing
binary variables and big-M parameters. This results in the following MILP form
of the master problem.
\begin{subequations}
\begin{align}
	& & \min & \; c^\top x + z \nonumber\\
        & &\text{s.t.} & \; z \geq \left(b^1 - A^1 x\right)^\top \eta^{1,t} +
        \left(b^2 - A^2 x\right)^\top\eta^{2,t}  + \bphiIP(b^2 - A^2 x) \eta^{\phi,t} +
        {\alpha}^t -  M_D(1 - u_t) \label{eq:lbf1}\\
        & & & \sum\limits_{t \in T} u_t = 1\label{eq:lbf2}\\
        & & &\ \bphiIP(b^2 - A^2 x) = \left(b^2 - A^2 x - G^2_I y^*_I\right)^\top{\eta^*} +
        {d^2_I}^\top y^*_I + M_P v \label{eq:ubf1}\\
        & & &\ -  \underline M v^1 \leq {\left(G_B^2\right)}^{-1} \left(b^2 -
        A^2 x - G^2_I y^*_I\right) \leq \bar M \left(1 - v^1\right) -
        \epsilon\label{eq:dr1}\\
        & & &\ M v - \sum\limits_{i \in \{1,\hdots,m_2\}} v^1_i \geq
        0\label{eq:dr2}\\
        & & &\ v - \sum\limits_{i \in \{1,\hdots,m_2\}} v^1_i \leq
        0\label{eq:dr3}\\
        & & &\ x \in X, z \text{ free}, \bphiIP \text{ free}, u_t \in \B \;
        \forall t \in T, v \in \B, v^1_i \in \B \; \forall i \in \{1, \hdots,
        m_2\}\nonumber
\end{align}
\end{subequations}

Constraints~\eqref{eq:lbf1}-\eqref{eq:lbf2} eliminate the minimization
operator by adding the binary variables $u_t$ for $t \in T$ and the big-M
parameter $M_D$. Constraints~\eqref{eq:ubf1}-\eqref{eq:dr3} eliminate the
if-else condition by adding the binary variables $v$, $v_i^1$ for $i \in \{1,
\hdots, m_2\}$ and the big-M parameters $\underline M$, $\bar M$ and $M$.
Specifically, the constraints~\eqref{eq:dr1}-\eqref{eq:dr3} impose domain
restriction on the first-stage variables with respect to the primal
function~\eqref{eq:milp-ubf}. This in turn restricts the domain of the
right-hand sides for which $\phiIP$ may be evaluated. The main idea is to set
some $v^1_i = 0$ whenever ${\left(G_{B_i}^2\right)}^{-1} \left(b^2 - A^2 x - G^2_I
y^*_I\right) \geq 0$ and to set $v^1_i = 1$ whenever ${\left(G_{B_i}^2\right)}^{-1}
\left(b^2 - A^2 x - G^2_I y^*_I\right) < 0$, where
${\left(G_{B_i}^2\right)}^{-1}$ is the $i^\textrm{th}$ row of the
optimal basis matrix inverse. If at least one $v^1_i = 1$, then $v = 1$
further implying $\bphiIP$ will have a very large value which is as
required. If all $v^1_i = 0$, then $v = 0$ further implying $\bphiIP$ will
have a finite value which is also as required. Finally, $\epsilon$ is also a
parameter corresponding to the domain restriction constraints added to deal
with the strict inequality arising from ``otherwise'' condition
in~\eqref{eq:milp-ubf}, and is the trickiest of all parameters to evaluate.
The discussion on finding appropriate values of big-M and $\epsilon$ parameters
is out of the scope of this paper.

\section{Conclusions \label{sec:conclusions}}
We have described a generalization of Benders' decomposition framework and
illustrated its principles by applying it to several well-known classes of
optimization problems that fall under the broad umbrella of MMILPs. The
development of an abstract framework for generalizing the principles of Benders'
technique for reformulation that encompasses non-traditional problem classes,
the specification of an associated algorithmic procedure, and its application to
the class of MIBLPs are our main contributions. These stemmed from our
observation that Benders' framework can be viewed as a procedure for iterative
refinement of dual functions associated with the value function arising from the
projection of the original problem into the space of first-stage variables, and
that this basic concept can be applied to a wide range of problems defined by
additively separable functions.

A conceptual extension of the generalized Benders' decomposition algorithm
from MIBLPs to the case of general MMILPs is straightforward. Similar
to MIBLPs, an $l$-stage MMILP can be formulated as a standard mathematical
optimization problem by considering a constraint requiring values of all but
first-stage variables to be optimal for an ($l-1$)-stage MMILP that is
parametric in the first-stage variables, in addition to the usual linear
constraints. Then, assuming that all input vectors and matrices are rational of
appropriate dimensions without loss of generality, we have an $l$-stage MMILP
with a parametric right-hand side $\beta$ defined as
\begin{equation}\label{eq:MMILPl}
\begin{aligned}
  \text{MMILP}^l(\beta) = \min & \; {d^{11}}^\top x^1 + {d^{12}}^\top x^2 +
  \hdots + {d^{1l}}^\top x^l\\
        \text{s.t.} & \text{ } A^{11}x^1 + A^{12} x^2 + \hdots + A^{1l} x^l \geq
        \beta \\
        & \; x^1 \in X^1 \\
        &\ (x^2, x^3, \hdots, x^l) \in \text{optimal set of
        $\text{MMILP}^{l-1}(b^2 - A^{21} x^1)$},
\end{aligned}
\end{equation}
where $\text{MMILP}^{l-1}(b^2 - A^{21} x^1)$ denotes an ($l-1$)-stage MMILP with
the parametric right-hand side $b^2 - A^{21} x^1$, which in turn is a linear
function of first-stage variables, defined similar to~\eqref{eq:MMILPl} but with
$l-1$ variable vectors, as
\begin{equation}\label{eq:MMILPl-1}
\begin{aligned}
  \text{MMILP}^{l-1}(\beta) = \min & \; {d^{22}}^\top x^2 + {d^{23}}^\top x^3 +
  \hdots + {d^{2l}}^\top x^l\\
        \text{s.t.} & \text{ } A^{22} x^2 + A^{23} x^3 + \hdots + A^{2l} x^l
        \geq \beta \\
        & \; x^2 \in X^2 \\
  &\ (x^3, x^4, \hdots, x^l) \in \text{optimal set of $\text{MMILP}^{l-2}(b^3 -
  A^{31} x^1 - A^{32} x^2)$}.
\end{aligned}
\end{equation}

These formulations exhibit the natural recursive property of MMILPs that we
spoke about in the beginning of the paper. It should be clear why this
recursive structure also means that the proposed framework makes it easy to
envision algorithms for solving such problems (whether these algorithms are
practical is another question). We can project~\eqref{eq:MMILPl} (for a fixed
$\beta = b^1$) into the space of the first-stage variables to obtain master and
subproblems. The subproblem itself involves an ($l-1$)-stage
MMILP~\eqref{eq:MMILPl-1}, and solution of it calls for solving this
($l-1$)-stage MMILP. This ($l-1$)-stage MMILP can as well be solved with the
generalized Benders' principle due to the recursive structure. Strong dual
functions can be constructed using techniques similar to those described
earlier in the paper.

Although not discussed here, this framework can be readily applied to even
broader classes of problems, such as those discussed in~\cite{BolConRalTah20},
which incorporate stochasticity. While it is unclear whether such algorithms
would be of practical interest, the algorithmic abstraction itself serves to
illustrate basic theoretical principles, such as concepts of general duality
and why $l$-stage MMILPs are canonical hard problems for stage $l$ of the
polynomial time hierarchy.

The algorithms described in this paper are naive in the sense that their
efficient implementations for practical purposes would require substantial
additional development, especially for the classes of MIBLPs and MMILPs. To
this end, our plans include enhancement of these algorithms by working in the
areas of preprocessing techniques, warm starting of master and subproblem
solves, cut management, a branch-and-Benders'-cut framework, alternative
linearization techniques, and other enhancements, as well as incorporating
aspects of these techniques into hybrid, non-Benders-type algorithms.

\section*{Acknowledgements}
This research was made possible with support from National Science
Foundation Grants CMMI-1435453, CMMI-0728011, and ACI-0102687, as well
as Office of Naval Research Grant N000141912330.

\bibliography{references}

\end{document}